%% file: main.tex
\documentclass[12pt,sort&compress]{elsarticle}
\input{mathhead}

\usepackage{hyperref}
\usepackage{fancyvrb}
\usepackage{fancyhdr}
\biboptions{sort&compress}

\usepackage{algorithm}
\usepackage[noend]{algpseudocode}
\newcounter{algsubstate}
\renewcommand{\thealgsubstate}{\alph{algsubstate}}

\usepackage{mathtools}
\mathtoolsset{showonlyrefs=true}

\usepackage{dsfont}
\usepackage{bm}
\usepackage{mathrsfs}
\usepackage{amsmath,mathtools}
\usepackage{amsfonts}
\usepackage{amsthm}
\usepackage{xcolor}
\usepackage{setspace}
\usepackage{float}
\usepackage{url}
\usepackage{multicol}
\usepackage{subfigure}
\usepackage{bm}
\usepackage[top=1in,bottom=1in,left=1in,right=1in]{geometry}

\definecolor{officegreen}{rgb}{0.0, 0.5, 0.0}

\newcommand{\mcS}{\mathcal{S}}

\newcommand{\mcV}{\mathcal{V}}
\newcommand{\mcL}{\mathcal{L}}

\newcommand{\mbRd}{{\mathbb{R}^d}}

\usepackage{multirow}

\begin{document}

\theoremstyle{definition}
\newtheorem{remark}{Remark}

\newcommand{\vertiii}[1]{{\left\vert\left\vert\left\vert #1
    \right\vert\right\vert\right\vert}}

\begin{frontmatter}

\title{Data-driven learning of robust nonlocal physics from high-fidelity synthetic data}

\address[yy]{Department of Mathematics, Lehigh University, Bethlehem, PA}
\address[nt]{Center for Computing Research, Sandia National Laboratories, Albuquerque, NM}
\address[md]{Computational Science and Analysis, Sandia National Laboratories, Livermore, CA}

\author[yy]{Huaiqian You}\ead{huy316@lehigh.edu}
\author[yy]{Yue~Yu\corref{cor1}}\ead{yuy214@lehigh.edu}
\author[nt]{Nathaniel Trask}\ead{natrask@sandia.gov}
\author[nt]{Mamikon Gulian}\ead{mgulian@sandia.gov}
\author[md]{Marta D'Elia}\ead{mdelia@sandia.gov}

\begin{abstract}

 A key challenge to nonlocal models is the analytical complexity of deriving them from first principles, and frequently their use is justified a posteriori. In this work we extract nonlocal models from data, circumventing these challenges and providing data-driven justification for the resulting model form. Extracting provably robust data-driven surrogates is a major challenge for machine learning (ML) approaches, due to nonlinearities and lack of convexity. Our scheme allows extraction of provably invertible nonlocal models whose kernels may be partially negative. To achieve this, based on established nonlocal theory, we embed in our algorithm sufficient conditions on the non-positive part of the kernel that guarantee well-posedness of the learnt operator. These conditions are imposed as inequality constraints and ensure that models are robust, even in small-data regimes. We demonstrate this workflow for a range of applications, including reproduction of manufactured nonlocal kernels; numerical homogenization of Darcy flow associated with a heterogeneous periodic microstructure; nonlocal approximation to high-order local transport phenomena; and approximation of globally supported fractional diffusion operators by truncated kernels.
\end{abstract}


\end{frontmatter}


\section{Background}

In contrast to partial differential equation (PDE) models which typically govern classical continuum mechanics and standard diffusion, nonlocal models describe systems in terms of integro-differential equations (IDEs). IDEs possess several modeling advantages over PDEs: reduced regularity requirements admit simpler descriptions of fracture mechanics\cite{ha2011characteristics,littlewood2010simulation,ouchi2015fully}; integral operators allow more natural description of long-range molecular interactions such as electrokinetic and surface tension effects \cite{fife2003some,tartakovsky2016smoothed,ansari2011calibration,seleson2009peridynamics,silling2000reformulation}; and compactly supported integral operators introduce a modeling lengthscale $\delta$ which may be used to model subgrid microstructures \cite{arambakam2014modeling,chebykin2011nonlocal,xu2008peridynamic,askari2008peridynamics,alali2012multiscale}. We consider in this work compactly supported nonlocal diffusion models, extensively analyzed in the literature \cite{javili2019peridynamics,du12sirev,Meerschaert2006}. While they have been effectively used, it is in general unclear how to justify the use of a kernel for a given application. As an example, in the peridynamic theory of continuum mechanics, the choice of nonlocal kernel is often justified a posteriori. Free parameters in constitutive models, including the lengthscale $\delta$, are often tuned to obtain agreement with experiment. In this manner, peridynamics has yielded enhanced descriptions of fracture mechanics, composite mechanics, but it is an open question how to obtain predictive nonlocal models from first principles. 

Recent work presents compelling justification of nonlocal kernels' role bridging microstructure and homogenized continuum response. Weckner and Silling characterize the multiscale response of elastic media by calibrating nonlocal dispersion relations to neutron scattering data \cite{weckner2011determination}. Delgoshaie argued that nonlocal diffusion kernels naturally describe tortuous flow pathways in porous media \cite{delgoshaie2015non}, while Chung obtained a nonlocal description of flow in fractured media via a multiscale finite element framework \cite{chung2018non}. Several authors have established conditions where nonlocal diffusion kernels admit interpretation as jump rates of stochastic processes  \cite{du2014nonlocal,defterli2015fractional,DElia2017}. Recently, Xu and Foster derived kernels describing an elastic medium with heterogeneous microstructure via homogenization\cite{xu2020deriving}. 

Works which provide rigorous a priori justification of nonlocal kernels typically require mathematically complex and lengthy derivations, often relying on restrictive assumptions of microstructure geometry. The goal of the current work is to establish a data-driven extraction of nonlocal models from observations of high-fidelity simulations, skirting the mathematical complexity of deriving nonlocal models. We pursue an inequality-constrained least squares approach which guarantees that extracted models are well-posed. In addition to learning unknown physics from systems, our method may be used to perform numerical homogenization, extracting efficient coarse-grained nonlocal models from high-fidelity synthetic data, such as microstructure-resolved local simulations. Furthermore, our algorithm may be applied to learn computationally cheap approximations of nonlocal operators of fractional-type while preserving accuracy. The latter are characterized by an infinite lengthscale and, as such, are computationally prohibitively expensive. Our technique delivers compactly supported nonlocal models that mimic the action of a fractional Laplacian at a much cheaper cost.

An additional open question in nonlocal modeling is whether kernels should be strictly positive. In some cases, nonlocal solutions can be interpreted as probability density functions of jump processes with jump rate given by the kernel, which, in turn, must be positive \cite{du2014nonlocal,defterli2015fractional,DElia2017}. On the other hand, in multiscale problems, several authors have extracted sign-changing kernels \cite{weckner2011determination}; see also the discussion in \cite{mengesha2013analysis}. It is unclear however, whether multiscale physics inherently lead to sign-changing kernels or if one could derive equally descriptive positive kernels. We note also that the mathematical analysis of non-positive kernels is more involved, suggesting that strictly positive kernels may be more desirable to work with \cite{mengesha2013analysis}. Recently, Mengesha and Du provided necessary conditions for diffusion processes governed by sign-changing kernels to be well-posed \cite{mengesha2013analysis}. Their analysis requires that the kernel may be decomposed into positive and non-positive parts, and that the negative part be sufficiently small so as not to compromise the coercivity of the overall operator. Some other recent works have explored extraction of nonlocal models from data, using either PDE-constrained optimization \cite{d2016identification,d2019imaging} or nonlocal generalizations of physics-informed neural networks \cite{pang2020npinns,Pang2019fPINNs} to obtain strictly positive kernels, which are thus well-posed. A unique feature of our approach is the extraction of an inequality constraint from Mengesha and Du's theory \cite{mengesha2013analysis} which allows learning of sign-changing kernels. Approximation of the kernel by Bernstein polynomials allows for a simple implementation of the constraint, guaranteeing extraction of well-posed sign-changing kernels. We will see that in various settings, sign-changing kernels sometimes appear naturally and sometimes do not, and therefore this data-driven derivation may suggest whether such kernels are appropriate for a given application.

\paragraph{Outline of the paper} The paper is organized as follows. To begin, in Section \ref{sec:abstract-regression} we define an abstract problem characterizing the fitting of nonlocal problems to high-fidelity observations. In Section \ref{sec:nonlocal-kernels} we review the well-posedness theory for sign-changing kernels that guides the discovery of well-posed nonlocal models and formulate the learning problem as an inequality-constrained optimization problem. Then, in Section \ref{sec:algorithm} we provide a two-stage algorithm to solve the optimization problem efficiently. In Section \ref{sec:manufactured_kernels} we illustrate properties of the proposed technique by using manufactured kernels; these results highlight consistency and robustness of our method. In Sections \ref{sec:darcy}, \ref{sec:homogenization}, and \ref{sec:fractional} we consider specific applications such as Darcy's flow, homogenization of a higher-order PDE and a fractional diffusion model. For each application, we highlight the impact of the proposed regression technique. Concluding remarks are reported in Section \ref{sec:conclusion}.

\section{Abstract problem}\label{sec:abstract-regression}

We consider in this section a framework for extraction of nonlocal operators which match observations of a given system in an optimal sense. We assume that a quantity of interest $u$ satisfies the problem
\begin{equation}\label{high-fidelity}
\left\{
\begin{aligned}
\mathcal{L}_{\text{HF}} [u](x) &= f(x) \quad x \in \Omega,\\
\mathcal{B}u(x)& = q(x) \quad x \in \partial \Omega, \text{ or } u \text{ is periodic,}
\end{aligned}
\right.
\end{equation}
where $\mathcal{L}_{\text{HF}}$ is a high-fidelity (HF) operator that faithfully represents the system; this can either be a partial differential equation or a fractional-differential operator. The function $f$ denotes a forcing term and operator $\mathcal{B}$ specifies a boundary or exterior condition. For example, when $\mathcal{L}_{\text{HF}}$ is a local operator, $\mathcal{B} = \text{Id}$ corresponds to a Dirichlet boundary condition and $\mathcal{B} = {\partial}/{\partial n}$ a Neumann boundary condition. We hypothesize that solutions to this HF problem may be approximated by solutions to an equivalent nonlocal volume-constrained problem of the form
\begin{equation}\label{coarsegrained}
\left\{
\begin{aligned}
\mathcal{L}_K[u](x) &= f(x) \quad x \in \Omega,\\
\mathcal{B}_{I} u(x) &= q(x) \quad  x \in \Omega_I, \text{ or } u \text{ is periodic.}
\end{aligned}
\right.
\end{equation}
Here, $\mathcal{L}_K$ is a nonlocal operator parametrized by $K$, defined below, $\Omega_I$ an appropriate nonlocal interaction domain, and $\mathcal{B}_I$ the corresponding nonlocal interaction operator specifying a volume constraint. 
In \eqref{coarsegrained}, $u$ and $f$ may coincide with the quantities of interest $u$ and $f$ in \eqref{high-fidelity}, or they may be appropriate ``coarse-grained'' representations of the same, we postpone a discussion of this multiscale scenario until Section \ref{sec:darcy}. 

We seek $\mathcal{L}_K$ as a linear nonlocal operator of the form
\begin{equation}\label{master}
\mathcal{L}_K[u](x) = -\int_{\Omega\cup\Omega_I} K(x,y) \left( u(y) - u(x) \right)\, dy,
\end{equation}
where $K$ denotes a nonlocal kernel compactly supported on the ball of radius $\delta$ centered at $x$, i.e. $B_\delta(x)$. Desired properties of the kernel and details on the analysis of problems associated with operator \eqref{master} are reported in Section \ref{sec:nonlocal-kernels}.
The interaction domain is defined as 
\begin{equation}\label{interaction-domain}
\Omega_I = \{ y \in \mathbb{R}^d \text{ such that } |y - x| \le \delta \text{ for some } x \in \Omega \};
\end{equation}
a condition for $u$ on $\Omega_I$ is required to specify $\mathcal{L}_K u(x)$ for $x \in \Omega$, by \eqref{master}. Thus, an appropriate exterior volume condition on $\Omega_I$ replaces the local boundary condition in \eqref{coarsegrained}. For simplicity we will assume $\mathcal{B}_I$ to be given: in the case of manufactured solutions one may apply the solution as a Dirichlet condition, while in the periodic case it need not be specified.

To learn the operator $\mathcal{L}_K$, we assume that we are given a collection of pairs of forcing terms $f_i$ and corresponding $u_i$ which arise from solutions to \eqref{high-fidelity},
\begin{equation}
\mathcal{D} = \left\{(u_i, f_i)\right\}_{i=1}^N.
\end{equation}
These observations will be used to train \eqref{master}.
In the simplest setting, $\mathcal{D}$ consists of pairs of forcing terms and corresponding solutions of \eqref{high-fidelity}. In specific coarse-graining applications, we can apply appropriate post-processing -- averaging, for example -- to coarsen such ``fine-scale'' solutions to construct $\mathcal{D}$. We then extract the nonlocal model by finding a kernel such that the action of $\mathcal{L}_K$ most closely maps $u_i$ to $f_i$.

That is, we solve an optimization problem of the form
\begin{equation}\label{abstractOptProblem}
K^* = \underset{K}{\text{argmin}} \sum_{i=1}^N \big\| \mathcal{L}_K [u_i] - f_i \big\|_{\mathcal{X}},
\end{equation}
where $\big\| \cdot \big\|_{\mathcal{X}}$ denotes an appropriate norm over $\Omega$.

In doing so, since $f_i$ represents $\mathcal{L}_\text{HF} [u_i]$ either exactly or in a post-processed sense, $\mathcal{L}_{K^*}$ best matches the action of $\mathcal{L}_{\text{HF}}$ on the training data $\mathcal{D}$, in the sense of problem \eqref{abstractOptProblem}.

In this manner, solution of the problem associated with $\mathcal{L}_{K^*}$ provides a surrogate for the given high-fidelity problem \eqref{high-fidelity}. We note, however, that there is no reason to expect the associated problem \eqref{coarsegrained} to be solvable. Previous works have imposed solvability via imposition of kernel positivity (either by positivity constraint \cite{d2016identification,d2019priori,d2019imaging} or restricting parameterization to only admit positive kernels \cite{pang2020npinns,Pang2019fPINNs}), they do not generalize naturally to kernels with negative part. In the next section, we consider a particular class of kernels allowing the addition of constraints to the optimization problem \eqref{abstractOptProblem} to guarantee extraction of well-posed models. 

We use this abstract setting as a framework for extracting nonlocal surrogates in a variety of scenarios where nonlocal models are expected to provide advantages. Before specializing to these application-specific settings, we first gather relevant theory regarding the stability and well-posedness of sign-changing nonlocal operators.

\section{Data-driven discovery of well-posed nonlocal models}\label{sec:nonlocal-kernels}
The ultimate goal of the operator regression discussed above is obtain a nonlocal model \eqref{coarsegrained} that can be solved for a class of forcing terms $f$ of interest -- for example, forcing functions at a suitable coarse scale and with certain regularity properties. However, without a careful definition and parametrization of a feasible set of such nonlocal operators $\mathcal{L}_K$ of the form \eqref{master}, there is no guarantee that regression will recover an operator $\mathcal{L}_K$ that enjoys this property, meaning that when \eqref{coarsegrained} is discretized, the resulting linear system will be singular.  Such operators cannot be used to extrapolate and generate predictions from new forcing functions $f$. Thus, we must avoid learning kernels $K^*$ that solve \eqref{abstractOptProblem}, but lead to non-invertible operators $\mathcal{L}_{K^*}$. 

The article \cite{mengesha2013analysis} gives an overview and establishes a generalized framework of well-posed nonlocal problems based on properties of the nonlocal kernel. In this context, well-posedness refers to existence and uniqueness of a solution to a nonlocal problem in an appropriate weak or variational sense, implying that suitable finite-dimensional discretizations of the problem can be solved numerically in a consistent and stable way. 

Below, we review the theory of \cite{mengesha2013analysis} which we will use throughout the article to ensure the discovery of well-posed nonlocal models. This reference discusses two sets of conditions for well-posedness. First, it reviews the case of nonlocal operators with \textit{nonnegative} kernels, reviewed in Section \ref{sec:nonnegative-kernels}. Next, it weakens this condition and establishes well-posedness for a class of sign-changing kernels satisfying conditions on the negative part of the kernel function, reviewed in Section \ref{sec:sign-changing-kernels}. We shall consider operator regression in the latter class of kernels. Our purpose for this is two-fold: first, this more general class of kernels is, in principle, more expressive and able to provide better fit to data. Second, we aim to study the utility of this  sign-changing kernels for nonlocal modeling. In other words, we study the question: is there an advantage to using sign-changing kernels, as opposed to nonnegative kernels, in fitting nonlocal models to data?

The conditions for well-posedness lead to constraints for the operator regression optimization problem. We define our parametrization of the unknown kernel in Bernstein polynomials and state the constrained optimization problem that arises in Section \ref{regression_algorithm}. 

\subsection{Well-posedness with nonnegative kernels}\label{sec:nonnegative-kernels}
In this section we introduce a standard nonlocal model with a nonnegative kernel and recall conditions for the well-posedness of associated diffusion problems. Let $\Omega\subset\mbRd$ be a connected bounded domain with sufficiently smooth boundary and let $\Omega_I$ be its interaction domain defined as in \eqref{interaction-domain}. We denote their union by $\widehat\Omega=\Omega\cup\Omega_I$. We consider the case of homogeneous nonlocal volume constraints\footnote{
While the authors of \cite{mengesha2013analysis} consider such volume constraints, their proofs suggest the analysis extends to the case of periodic conditions. We consider both homogeneous Dirichlet and periodic conditions in our examples in Sections \ref{sec:manufactured_kernels} -- \ref{sec:fractional}.} on $\Omega_I$, i.e. $\mathcal B_I u=0$. We define the action of a nonlocal diffusion operator $\mathcal{L}_\rho$ on a function $u$ as
\begin{equation}\label{L}
\mcL_\rho u (x) = -2 \int_{\widehat\Omega} \rho(|y-x|)
\left( u(y) - u(x) \right) dy,    
\end{equation}
where the nonnegative and compactly supported kernel $\rho= \rho(|\xi|)$ is such that
\begin{equation}\label{rho-conditions}
\begin{array}{l}
|\xi|^2 \rho(|\xi|) \in L^1_{\text{loc}}(\mathbb{R}^d), \text{ and}\\[2mm]
\exists\, \sigma > 0 \text{ such that } (0,\sigma) \subset \text{supp}(\rho).
\end{array}
\end{equation}
To study existence and uniqueness of solutions to the equation
\begin{equation}\label{eq:nonnegative_strong_form}
\mathcal{L}_\rho u(x) = b(x), \quad {\rm for} \;\; x\in\Omega,
\end{equation}
we define the space $S({\widehat\Omega})$ associated with $\mcL_\rho$ as the typical energy space considered for nonlocal problems \cite{gule10,du12sirev,du13m3as,mengesha2013analysis}, i.e.
\begin{displaymath}
S({\widehat\Omega}) = \left\{u \in L^2({\widehat\Omega}) : 
\int_{\widehat\Omega} \int_{\widehat\Omega} \rho(|y-x|) |u(y)-u(x)|^2 dy dx < \infty
\right\}.
\end{displaymath}
Note that $S({\widehat\Omega})$ arises naturally from the symmetric bilinear form
\begin{align}
(\mathcal{L}_\rho u, w)_{L^2({\Omega})} 
&=
-2 \int_{{\Omega}} w(x) \int_{{\widehat\Omega}} \rho(|y-x|) \left(u(y)-u(x)\right) dy dx \\
&= 
\int_{\widehat\Omega} \int_{\widehat\Omega} \rho(|y-x|) \left(u(y)-u(x)\right) \left(w(y)-w(x)\right) dy dx,
\end{align}
where $w$ is such that $w=0$ in $\Omega_I$. We denote by $\mcS'$ the dual space of $\mcS$. Let $V$ be a closed subspace of $L^2({\widehat\Omega})$ that contains the constant function $u \equiv 0$ and no other constant functions, we define the solution space as $\mcV_s=V\cap\mcS(\widehat\Omega)$. According to Corollary 1 of \cite{mengesha2013analysis}, if $\rho$ satisfies \eqref{rho-conditions}, there exists a coercivity constant $\kappa = \kappa(\rho, V, {\widehat\Omega})$ of $\mathcal{L}_\rho:\mcS\to\mcS'$ such that
\begin{equation}\label{L-coercivity}
\| u \|_{L^2({\Omega})}^2 \le \kappa \left( \mathcal{L}_\rho u, u \right)_{L^2({\Omega})}, \quad
\text{ for all } u \in V \cap S({\widehat\Omega}),
\end{equation}
where $(\cdot,\cdot)_{L^2(\cdot)}$ denotes the $L^2$ inner product. By Lemma 3.1 of \cite{mengesha2013analysis}, this result guarantees that there exists a unique $u\in\mcV_s$ such that, for $b\in\mcS'$ and for all $w\in \mcV_s$
\begin{equation}\label{positive-weak-inner-product}
(\mathcal{L}_\rho u , w )_{L^2(\Omega)} = ( b , w )_{L^2(\Omega)},
\end{equation}

Equation \eqref{positive-weak-inner-product} is the weak form of \eqref{eq:nonnegative_strong_form}. 

\subsection{Well-posedness for sign-changing kernels}\label{sec:sign-changing-kernels}
In this section we consider the type of sign-changing kernels studied by \cite{mengesha2013analysis}. These kernels define operators of the form
\begin{equation}\label{Plambda}
\mathcal{L}_K = \mathcal{L}_\rho + \lambda \mathcal{L}_g,
\end{equation}
where $\mcL_\rho$ is the operator defined in \eqref{L} and the nonlocal operator $\mathcal{L}_g$ on the right-hand side is defined as
\begin{align}
\mathcal{L}_g u (x) &= -2 \int_{\widehat\Omega} g(|y-x|)
\left( u(y) - u(x) \right) dy.
\end{align}
Thus, we consider integral operators of the form
\begin{equation}
\mathcal{L}_K u(x) =
- \int_{\widehat\Omega} K(x,y)
\left( u(y) - u(x) \right) dy, 
\end{equation}
with kernel
\begin{equation}\label{kernel}
K(x,y) = 2 \rho(|y-x|) + 2 \lambda g(|y-x|). 
\end{equation}

The function $\rho$ is assumed to satisfy the conditions \eqref{rho-conditions} for well-posedness of the problem \eqref{eq:nonnegative_strong_form}, so that the associated bilinear form is coercive as in \eqref{L-coercivity} with coercivity constant $\kappa$. 
The idea of \cite{mengesha2013analysis} is to think of  $\mathcal{L}_K$ as a small perturbation of $\mathcal{L}_\rho$ by $\lambda \mathcal{L}_g$. Then, for well-posedness of equations of the form 
\begin{equation}\label{eq:strong_form}
\mathcal{L}_K u = b,
\end{equation}
we need the contribution of the perturbation $\mathcal{L}_g$ not to compromise the coercivity of the variational form associated with $\mathcal{L}_\rho$, so that $\mathcal{L}_K$ will also induce a coercive, hence well-posed, variational problem. 

We now summarize the conditions on $g$ and $\lambda$ that guarantee well-posedness of the perturbed system \eqref{Plambda}. From Lemma 4.2, Theorem 4.3, and Corollary 2 of \cite{mengesha2013analysis}, 
have that if $g$ is a compactly supported locally integrable radial function, and if $g$ and $\lambda$ satisfy
\begin{equation}\label{eq:g-condition}
|\lambda| < \frac{1}
{2 \kappa \left( \|g\|_{L^1(\mathbb{R}^d)} + \|G\|_{L^\infty(\widehat\Omega)}\right)},
\end{equation}
with $G(x) = \int_{\widehat\Omega} g(|y-x|) dy$ 
and $\kappa$ defined as in \eqref{L-coercivity}, then, for $b \in \mcV_s'$, the weak form of equation \eqref{eq:strong_form}
\begin{equation}\label{sign-changing-weak}
(\mathcal{L}_K u , w )_{L^2(\Omega)} = ( b , w )_{L^2(\Omega)}, \quad \forall\,w \in \mathcal{V}_s,
\end{equation}
has a unique solution $u\in \mcV_s$. This is the analogue of \eqref{positive-weak-inner-product} for the equation \eqref{eq:strong_form}.

In the following section, we will parametrize $\rho$ and $\lambda$ in coefficients of a polynomial expansion, making $\lambda$ in \eqref{kernel}  a redundant parameter that scales the coefficients. Writing $h = \lambda g$, multiplying \eqref{eq:g-condition} through by the factor $\left( \|g\|_{L^1(\mathbb{R}^d)} + \|G\|_{L^\infty(\widehat\Omega)}\right)$ and pulling $\lambda$ inside the norms, the setup of \eqref{kernel} with condition \eqref{eq:g-condition} can be written as
\begin{align}\label{final_k_setup}
\begin{cases}
& K(x,y) = 2\rho(|y-x|) + 2h(|y-x|), \\
& H(x) = \int_{\widehat\Omega} h(|y-x|) dy, \\
& \|h\|_{L^1(\mathbb{R}^d)} + \|H\|_{L^\infty(\widehat\Omega)} < \frac{1}{2 \kappa }.
\end{cases}
\end{align}
Accordingly, we also define $\mathcal{L}_h = \lambda \mathcal{L}_g$, so that the operator $\mathcal L_K$ in \eqref{Plambda} can be written as $\mathcal L_K = \mathcal{L}_\rho + \mathcal{L}_h$. 

\subsection{Sign-changing kernel regression}\label{regression_algorithm}

\medskip
With the purpose of applying the abstract formulation introduced in Section \ref{sec:abstract-regression} to the discovery of kernels described in Section \ref{sec:nonlocal-kernels}, we introduce a representation of the kernel as a linear combination of Bernstein basis polynomials. Notably, we assume that $K$ is radial, i.e., that $K(x,y)$ is a function of $|\frac{x-y}{\delta}|$, and supported on $B_\delta(0)$, the ball centered at $0$ of radius $\delta$. On this support, we expand in polynomials in the radial variable. Specifically, we parametrize \eqref{kernel} by defining
\begin{equation}\label{model-kernel}
K(x,y) =
    \sum_{m=0}^{M}\frac{C_m}{\delta^{d+2}} B_{m,M}
    \bigg (\bigg|\frac{x-y}{\delta}\bigg|\bigg)
    +
    \sum_{m=0}^{M}\frac{D_m}{\delta^{d+2}} B_{m,M}
    \bigg (\bigg|\frac{x-y}{\delta}\bigg|\bigg),
\end{equation}
where the Bernstein basis functions are defined as
\begin{equation}
B_{m,M}(x) = \begin{pmatrix}
M\\
m\\
\end{pmatrix}
x^m(1-x)^{M-m}.
\end{equation}
This corresponds to setting
\begin{equation}\label{eq:rho-g-def}
2\rho(|\xi|)
=
\sum_{m=0}^{M}\frac{C_m}{\delta^{d+2}} 
B_{m,M}\bigg (\bigg|\frac{\xi}{\delta}\bigg|\bigg), \quad{\rm and} \quad
2h(|\xi|)
=
\sum_{m=0}^{M}\frac{D_m}{\delta^{d+2}} 
B_{m,M}\bigg (\bigg|\frac{\xi}{\delta}\bigg|\bigg)
\end{equation}
in \eqref{final_k_setup},
where the coefficients of the linear combinations are unknown. The scaling of $C_m$'s and $D_m$'s is standard and has been used in the literature (see, e.g. \cite{Chen2011}) to guarantee convergence of nonlocal diffusion operators to $\Delta$, the classical Laplacian, as $\delta\to 0$.

For $\mathbf{C} \in \mathbb{R}^{M+1}$, $\mathbf{D} \in \mathbb{R}^{M+1}$, the kernel discovery problem can be stated as
\begin{equation}\label{eq:true_problem}
\left\{
\begin{aligned}
&\underset{\mathbb{R}^{2M+2}}{\text{min }} \ L\left(\begin{bmatrix} \mathbf{C} \\ \mathbf{D} \end{bmatrix}\right) \\
&\text{subject to:}\quad \mathbf{C} \ge \mathbf{0},\quad N(\mathbf{D}) \le \left[ 2 \kappa(\mathbf{C}) \right]^{-1},
\end{aligned}
\right.
\end{equation}
where
\begin{equation}\label{eq:true_constraints}
\left\{
\begin{aligned}
&L\left(\begin{bmatrix} \mathbf{C} \\ \mathbf{D} \end{bmatrix}\right) = \frac{1}{N} \sum_{i=1}^N  \|(  \mathcal{L}^{\mathbf{C}} + \mathcal{L}^{\mathbf{D}}) [u_i] - f_i \|_{\mathcal{X}} ^2\\
&N(\mathbf{D}) = 
\frac{1}{2}
\left\|
\sum_{m=0}^M\frac{D_m}{\delta^{d+2}} 
B_{m,M}\bigg (\bigg|\frac{x}{\delta}\bigg|\bigg)
\right\|_{L^1(\mathbb{R}^d)}
+
\frac{1}{2}
\left\|
\sum_{m=0}^M\frac{D_m}{\delta^{d+2}} 
\int_{\mathbb{R}} B_{m,M}\bigg (\bigg|\frac{y-x}{\delta}\bigg|\bigg) dy
\right\|_{L^\infty(\widehat\Omega)}\\
&\text{$\kappa(\mathbf{C})$ is a constant satisfying \eqref{L-coercivity}.}
\end{aligned}
\right.
\end{equation}
Above, $\mathcal{L}^{\mathbf{C}}$ denotes $\mathcal{L}_{\rho}$,  $\mathcal{L}^{\mathbf{D}}$ denotes $\mathcal{L}_{h}$, following the parametrization \eqref{eq:rho-g-def}, and we delay specification of the norm $\|\cdot\|_{\mathcal{X}}$ until the following section. Numerical optimization of this problem is discussed in Section \ref{sec:algorithm}.

\section{Algorithm}\label{sec:algorithm}
The problem \eqref{eq:true_problem} is expected to be nonconvex and susceptible to many local minima. At the same time, while the constraint 
$\mathbf{C} \ge \mathbf{0}$ is easier to enforce, the second constraint involving $N(\mathbf{D})$ and $\kappa(\mathbf{C})$ in \eqref{eq:true_constraints} is highly nonlinear in both the left-hand and right-hand sides, posing a significant challenge. For numerical efficiency, we have followed a ``two-stage'' strategy for numerical optimization. 

We first find a set of $\mathbf{C^*}$ representing an initial fit to the data using only the nonnegative kernel $\rho$ as in \eqref{eq:rho-g-def}, ignoring the term $h$ and avoiding the nonlinear constraint. We minimize the loss
\begin{align}\label{eq:L1}
L_1(\mathbf{C})
&=
\frac{1}{N}\sum\limits_{i=1}^{N} \left\|\sum_{m=0}^M\frac{\text{ReLU}(C_m)}{\delta^3} \Bigg[\int_{\widehat{\Omega}}B_{m,M}\bigg (\bigg|\frac{y-x}{\delta}\bigg|\bigg)(u_i(x)-u_i(y))dy\Bigg]_{\!\!\Delta}-f_i(x)\right\|^2_{\mathcal{X}}.
\end{align}

Here we use the notation $\left[ \int \cdot \; dy\right]_\Delta$ to denote the discrete approximation of a nonlocal kernel by applying quadrature. For the purposes of this work, we will use Silling's one-point quadrature \cite{silling_2005_2}, popular in particle discretizations of peridynamics. We note however that any choice of quadrature for discretization of the integral in the strong form \cite{trask2019asymptotically,you2019asymptotically}, or basis for the discretization of the weak form, may be applied here; see \cite{ACTAnumerica2020} and references therein. For a given nonlocal operator $\mathcal{L}_K$ and basis functions $\mathbf{\Phi} = \left[\phi_p\right]_{p=1}^P$ and corresponding degrees of freedom $\mathbf{w} = \left[w_p\right]_{p=1}^P$ the discretization provides mass and nonlocal stiffness matrices 
\begin{equation}\label{eq:discretization-matrices}
M_{ij}= \left[ \int_\Omega \phi_i(x) \phi_j(x)\,dx \right]_\Delta
\quad {\rm and} \quad
S_{ij}= \left[\int_\Omega \mathcal{L}_K [\phi_i(x)] \phi_j(x) dx\right]_\Delta \!\!\!.
\end{equation} 

Note the presence of the $\text{ReLU}$ function in $L_1$. We minimize $L_1$ using the Adam optimizer \cite{kingma2014adam}, also mapping 
\begin{equation}
\mathbf{C} \mapsto \text{ReLU}(\mathbf{C})
\end{equation}
after each step of gradient descent. This leaves $L_1$ invariant but ensures that the sequence and local minimum produced by gradient descent satisfies $\mathbf{C} \ge 0$. Thus, the simpler constraint $\mathbf{C} \ge 0$ is hard-coded into the loss function and algorithm, allowing us to obtain $\mathbf{C}^*$ by applying the Adam algorithm to the unconstrained optimization  problem. 

We then compute the coercivity constant $\kappa(\mathbf{C}^*)$ corresponding to the operator $\mathcal{L}^{\mathbf{C}^*}$. For $\mathbf{w},S,$ and $M$ consistent with the choice of discretization $\left[ \int \cdot \; dy\right]_\Delta$, we solve the following generalized eigenvalue problem
\begin{equation}\label{eq:generalized_eigenvalue}
    \left[ \kappa(\mathbf{C}^*) \right]^{-1} = \min_{\mathbf{w}} \frac{\mathbf{w}^{T} S \mathbf{w}}{\mathbf{w}^{T}M\mathbf{w}}.
\end{equation}
This computation happens only \textit{once} in our algorithm, so that $\mathbf{C}^*$ and $\kappa(\mathbf{C}^*)$ are fixed in the next stage. 

In the second stage, we correct the initial fit given by $\mathcal{L}^{\mathbf{C}^*}$ by finding $\mathcal{L}^{\mathbf{D}^*}$ subject to the second constraint in \eqref{eq:true_problem}, in which the right-hand side $\left[ 2 \kappa(\mathbf{C}^*) \right]^{-1}$ of the constraint is fixed. We apply the augmented Lagrangian method \cite{yu2019dag,nocedal2006numerical} to solve for $\mathbf{D}$ under the constraint condition, using the function
\begin{align}
H(\mathbf{D},\theta)
&= 
\dfrac{1}{2\kappa(\mathbf{C}^*)}
-
N(\mathbf{D})
-\theta^2 \\
&\begin{multlined}[t]=
\dfrac{1}{2\kappa(\mathbf{C}^*)}
-
\frac{1}{2}
\sum_{m=0}^M\frac{D_m}{\delta^{d+2}} 
\left[
\int_{\mathbb{R}^d}B_{m,M}\bigg (\bigg|\frac{x}{\delta}\bigg|\bigg)
dy\right]_{\Delta}\\
-
\frac{1}{2}\sup_{\widehat{\Omega}}
\left|
\sum_{m=0}^M\frac{D_m}{\delta^{d+2}} 
\left[\int_{\mathbb{R}^d} B_{m,M}\bigg (\bigg|\frac{y-x}{\delta}\bigg|\bigg) dy\right]_{\Delta}
\right|
-\theta^2.
\end{multlined}
\end{align}
Here, $\theta$ is the slack variable arising from the inequality constraint. 
We then apply the Adam optimizer to the penalized loss function
\begin{align}\label{eq:L2}
L_2(\mathbf{D}, \theta)
&=
\dfrac{1}{N}\sum\limits_{i=1}^{N} 
\|
\mathcal{L}^{\mathbf{C}^* + \mathbf{D}} [u_i] -f_i \|^2_{\mathcal{X}}
+
\lambda H(\mathbf{D},\theta)+\dfrac{\mu}{2} H^2(\mathbf{D},\theta)
\\
&=
\begin{multlined}[t]
\frac{1}{N}\sum\limits_{i=1}^{N} \left\|\sum_{m=0}^M\frac{C_k + D_k}{\delta^3} \Bigg[\int_{\widehat{\Omega}}B_{m,M}\bigg (\bigg|\frac{y-x}{\delta}\bigg|\bigg)(u_i(x)-u_i(y))dy\Bigg]_{\!\!\Delta}-f_i(x)\right\|^2_{\mathcal X}\\
+
\lambda H(\mathbf{D},\theta)+\dfrac{\mu}{2}H^2(\mathbf{D},\theta),
\end{multlined}
\end{align}
adjusting the parameters $\lambda$ and $\mu$ as described in Algorithm \ref{alg:augmented_lagrangian}. We denote the minimum by $\mathbf{D}^*$.

Although this two-stage optimization algorithm is numerically efficient, there is a degree of overconstraining arising from a single computation of the coercivity constant $\kappa$. That is, although $\mathbf{C}^*$ is a local minimum of $L_1$ and $\mathbf{D}^*$ is a local minimum of $L_2$, $(\mathbf{C}^*, \mathbf{D}^*)$ is not necessarily a local minimum of $L$. 

A more advanced algorithm could involve iterated computation of the coercivity constraint $\kappa$. In our examples below, however, we found that the two-stage algorithm is sufficient and that there was no significant benefit to more complex and expensive optimization algorithms. Therefore, we have used the two-stage algorithm throughout. 
 
\begin{algorithm}\label{alg:regression}
\caption{Nonlocal kernel regression}\label{alg:augmented_lagrangian}
\begin{algorithmic}[1]
\State Initialize $C_m^{(0)} \sim \mathcal{U}\left(0,1\right)$.
\State
Obtain $\mathbf{C}^*$ as a local minimum of $L_1(\mathbf{C})$,
using the Adam optimizer while
updating $\mathbf{C} \leftarrow \text{ReLU}(\mathbf{C})$ after each step of gradient descent. 
\State
Select basis functions $\mathbf{\Phi}$ and corresponding degrees of freedom  $\mathbf{w}$, 
and assemble the matrices $M$ and $S$ in \eqref{eq:discretization-matrices}. 
\State
Solve the generalized eigenvalue problem \eqref{eq:generalized_eigenvalue} for $\kappa(\mathbf{C}^*)$.

\State
Initialize ${D}_m^{(0)} \sim \mathcal{U}\left(-\frac{1}{\sqrt{M+1}}, \frac{1}{\sqrt{M+1}}\right)$ and $\theta^{(0)}=1$.
\State Set $\text{\texttt{STEP\_{MAX}}}=100$, $\lambda=0$, $\mu=1$, $s=1$, $\rho=10$, $c=1/4$, $\epsilon=10^{-8}$.

\While{$s \leq \text{\texttt{STEP\_{MAX}}}$:} \Comment{Perform Augmented Lagrangian Algorithm for $\mathbf{D}$}
\State\label{step_solve_optimization}
Solve the unconstrained optimization problem
\begin{align*}
    &(\mathbf{D}^{(s)},\theta^{(s)})=\underset{\mathbf{D},\theta}{\text{argmin }} 
    L_2(\mathbf{D}, \theta).
\end{align*}
\If{$H(\mathbf{D}^{(s)},\theta^{(s)})\leq \epsilon$}
\State
Stop.
\Else
\If{$H(\mathbf{D}^{(s)},\theta^{(s)})\geq c H(\mathbf{D}^{(s-1)},\theta^{(s-1)})$}
\State
Update penalty $\mu\leftarrow \rho\mu$.
\If{$\mu\geq 10^{20}$}
\State
Stop.
\EndIf
\Else
\State
Update Lagrange multiplier $\lambda \leftarrow \lambda+\mu H(\mathbf{D}^{(s)},\theta^{(s)})$. 
\EndIf
\EndIf
\State
Update the iteration number $s\leftarrow s+1$.
\EndWhile
\State
$\mathbf{D}^* = \mathbf{D}^{(s)}$.
\end{algorithmic}
\end{algorithm}
In applying Algorithm \ref{alg:augmented_lagrangian}, we run the Adam optimizer (in PyTorch) using a batch size of 100 and learning rate 5e-3 throughout the article. We run until the loss stagnates, indicating that a stationary point has been reached. This was typically between 200 and 500 epochs for the first stage of the algorithm (to find $\mathbf{C}^*$) and 10 epochs for the second stage (to find $\mathbf{D}^*$) per iteration of the augmented Lagrangian method. In general we will take at least $O(1,000)$ samples for the purposes of training, to ensure we are well into the regime of having sufficient data. Throughout, we use the norm
\begin{equation}\label{eq:X_norm}
\| f \|_{\mathcal{X}} = 
\sqrt{ \frac{1}{\#\mathcal{X}} 
\sum_{x_k \in \mathcal{X}} f(x_k)^2 }
\end{equation}
in \eqref{eq:L1} and \eqref{eq:L2},
for a discrete collection $\mathcal{X}$ of points in $\Omega$. We specify the collection $\mathcal{X}$ in the sections below. 

\section{Tests with manufactured kernels}\label{sec:manufactured_kernels}
In this section, we give the first illustration of the method described above by training on data $\mathcal{D}_{\text{train}} = \{(u_i,f_i)\}_{i=1}^N$ generated from a nonlocal equation of the form \eqref{coarsegrained} with periodic boundary conditions, for a manufactured kernel $K_{\text{man}}$ that is given \textit{a priori}. This allows us to validate the approach and quantify its ability to handle nonnegative and sign-changing kernels, before moving to more realistic data in subsequent sections.

We generate 50,000 training pairs $(u_i,f_i)$ in the following way. First we randomly generate the coefficients of a Fourier series for $u_i$ as 
\begin{equation}\label{eq:low_freq_distribution}
\widehat{u}_k \sim \exp(-\alpha k^2) \xi, \quad \xi \sim \mathcal{U}[0,1],
\end{equation}
where $\mathcal{U}[0,1]$ denotes the uniform distribution on $[0,1]$ and $\alpha = 0.1$.
Then, $u_i$ is given by
\begin{equation}\label{eq:f_manufactured}
    u_i =
    \sum_{k=0}^{100} \widehat{u}_k \cos(2\pi kx/L)
\end{equation}
and $f_i$ is obtained from \eqref{high-fidelity} via a numerical Fourier transform.
We train the kernel using $\mathcal{X}$ in \eqref{eq:X_norm} to be 101 equidistant points in $[0,1]$.

The learned kernel $K^*$ is validated on a test set consisting of a different set of samples $\mathcal{D}_{\text{test}} = \{(u_j,f_j)\}_{j=1}^N$ generated as above, with $N=10,000$. We test two types of basis expansions in \eqref{model-kernel} for learning $K$: a linear kernel ($M=1$) and a quadratic kernel ($M=2$).

\begin{figure}[htpb!]
    \centering
    \includegraphics[width = 0.85\textwidth]{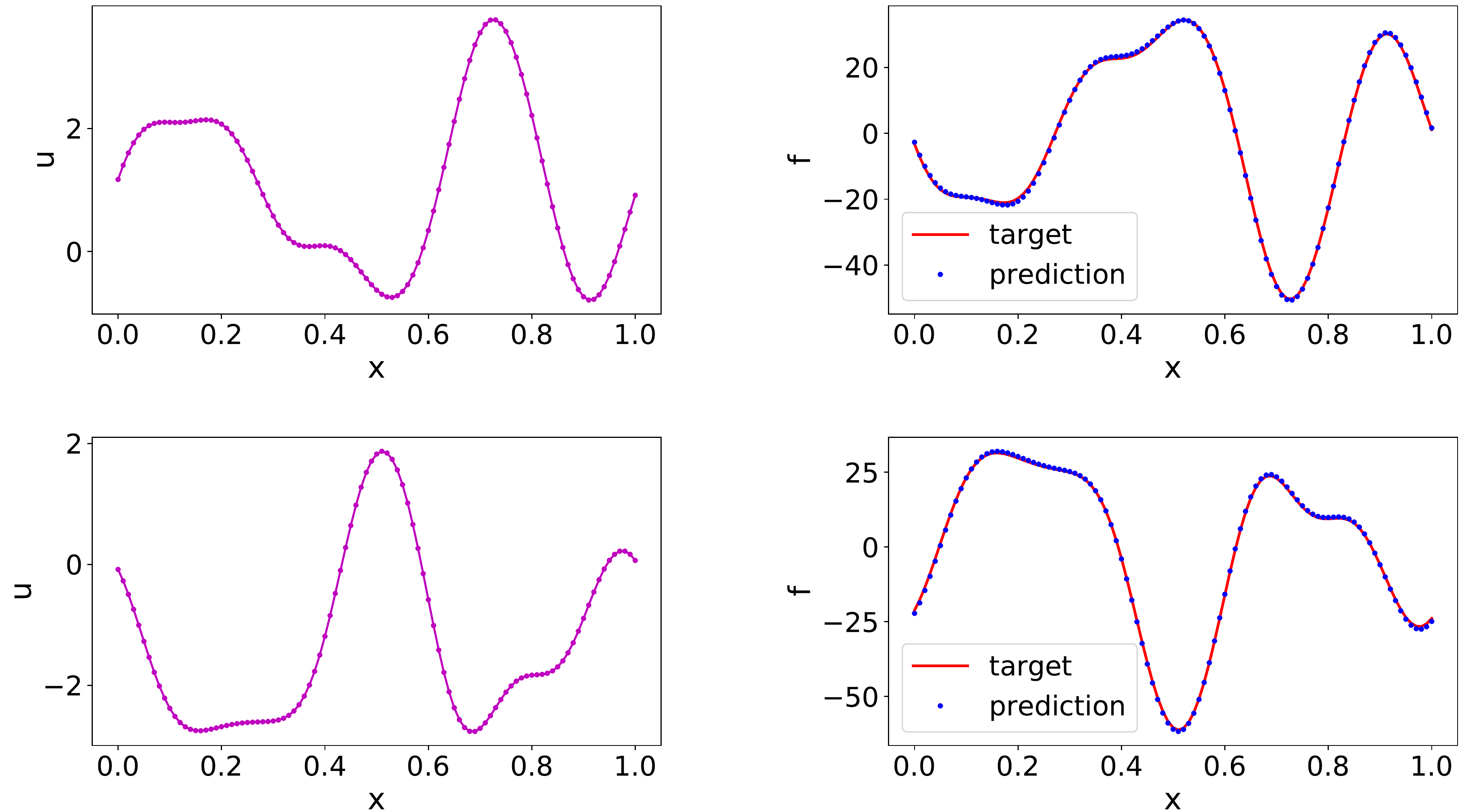}
    \caption{Comparison of the prediction $\mathcal{L}_{K^*}[u]$ and target $f$ when generating data from the manufactured kernel \eqref{eq:linearkernel} and fitting a nonnegative linear basis polynomial kernel $K^*$. The left column illustrates $u$, whereas the right column $f$ (target) and $\mathcal{L}_{K^*}[u]$ (prediction), for a sample of $(u,f)\in\mathcal D_{\rm test}$ in each row.}
    \label{fig:linear_linearbase}
\end{figure}

\begin{figure}[htpb!]
    \centering
    \includegraphics[width =0.85\textwidth]{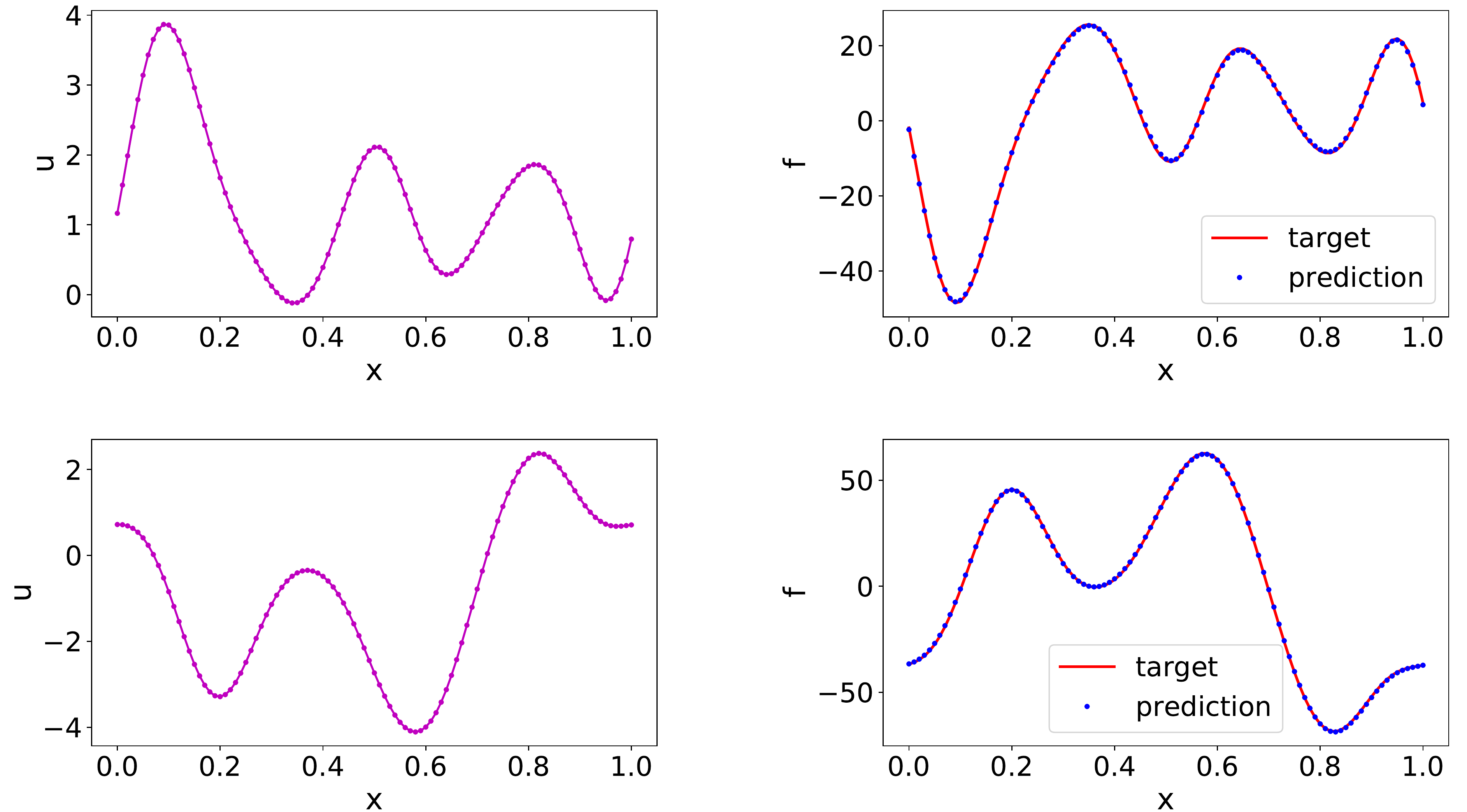}
    \caption{Comparison of the prediction $\mathcal{L}_{K^*}[u]$ and target $f$ when generating data from the manufactured kernel \eqref{eq:linearkernel} and fitting a nonnegative quadratic basis polynomial kernel $K^*$. The left column illustrates $u$, whereas the right column $f$ (target) and $\mathcal{L}_{K^*}[u]$ (prediction), for a sample of $(u,f)\in\mathcal D_{\rm test}$ in each row.}
    \label{fig:linear_quabase}
\end{figure}

For the case of a linear manufactured kernel
\begin{equation}\label{eq:linearkernel}
K_{\text{man}}(x,y) = \frac{4}{\delta^3}\bigg|\frac{y-x}{\delta}\bigg| \; \mathds{1}_{[0,1]}
    \Bigg(\bigg|\frac{y-x}{\delta}\bigg|\Bigg),
\end{equation}
used to generate data, prediction using $\mathcal{L}_{K^*}$ is shown for two test pairs $(u,f)$ in Figure \ref{fig:linear_linearbase} for regression using a linear kernel, and Figure \ref{fig:linear_quabase} for regression using a quadratic basis.  
In these first examples, we only learned the nonnegative coefficients $\mathbf{C}$, effectively setting $\mathbf{D}=0$ and not performing the second stage of the algorithm. This verifies that given data from a nonnegative manufactured kernel, the algorithm can provide a good fit with a nonnegative learned kernel. 

In the next example, we fit data generated from the sign-changing kernel
\begin{equation}\label{eq:posnegkernel}
    K_{\text{man}}(x,y) = \frac{21.4615}{\delta^3} \cos\bigg(\frac{3\pi|y-x|}{5\delta}\bigg) \; 
    \mathds{1}_{[0,1]}
    \Bigg(\bigg|\frac{y-x}{\delta}\bigg|\Bigg).
\end{equation}
Here, we apply the second stage of Algorithm \ref{alg:augmented_lagrangian}, learning both $\mathbf{C}$ and $\mathbf{D}$.  
We use this example to explore whether learning a kernel $K^*$ by fitting the action of $\mathcal{L}_K$ on the data $\mathcal{D}_{\text{train}}$ can be expected to reproduce the kernel $K_{\text{man}}$.
Of course, for a general kernel $K_{\text{man}}$, this can only be expected for increasing basis order $M$ of the polynomial expansion in \eqref{model-kernel}.
We generate two sets of training and test data. The first set, which we refer to as  ``low-frequency", consists of 50,000 samples of the form \eqref{eq:f_manufactured} with distribution \eqref{eq:low_freq_distribution}, just as for the previous example. The second set, which we refer to as ``high-frequency'', consists of $25,000$ of such low-frequency samples and $25,000$ samples $(u_i, f_i)$ where
\begin{equation}\label{eq:high_freq_distribution}
    u_i = \xi_1\sin(2\pi k_1 x/L)+
    \xi_2\cos(2\pi k_2 x/L),
\end{equation}
for $\xi_1, \xi_2 \sim \mathcal{U}[0,1]$ and $k_1$ and $k_2$ being random integers sampled uniformly in $[5,15]$, and with $f_i$ computed using a numerical Fourier transform. 
We learn kernels $K$ of the form \eqref{model-kernel} with degrees $M$ of the basis increasing $2$ to $20$, for both data sets. 

\begin{figure}[H]
    \centering
    \includegraphics[width=0.48\textwidth]{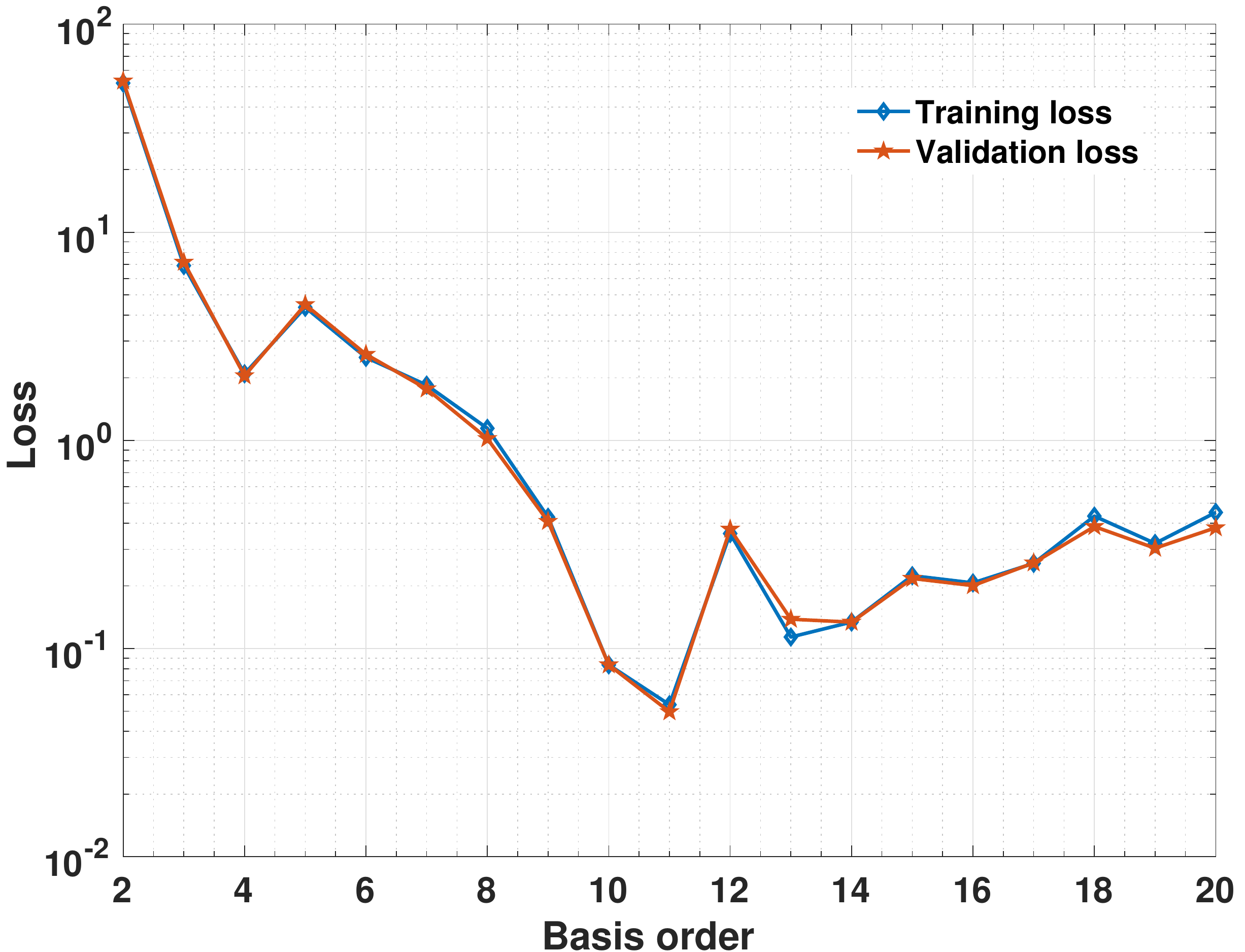}
    \includegraphics[width=0.48\textwidth]{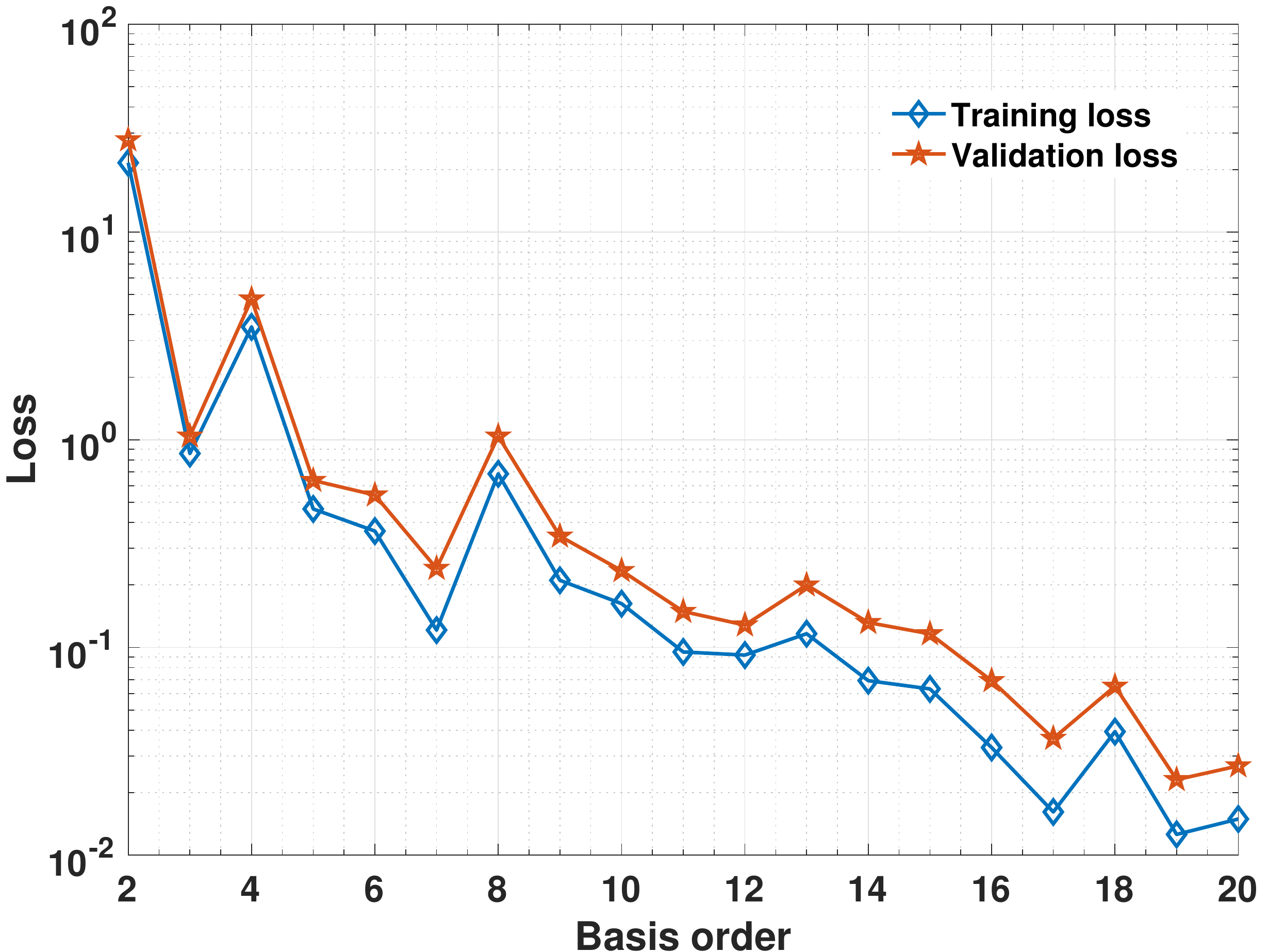}
    \caption{\textit{Left:} The training loss and validation loss versus basis order for low-frequency data generated using the manufactured kernel $K_{\text{man}}$ \eqref{eq:posnegkernel}.
    \textit{Right:} The training loss and validation loss versus basis order for high-frequency data generated using the manufactured kernel $K_{\text{man}}$ \eqref{eq:posnegkernel}.}
    \label{fig:myrelulossvsorder_and_Loss_Hfreq}
\end{figure}

\begin{figure}[H]
    \centering
    \includegraphics[width=0.48\textwidth]{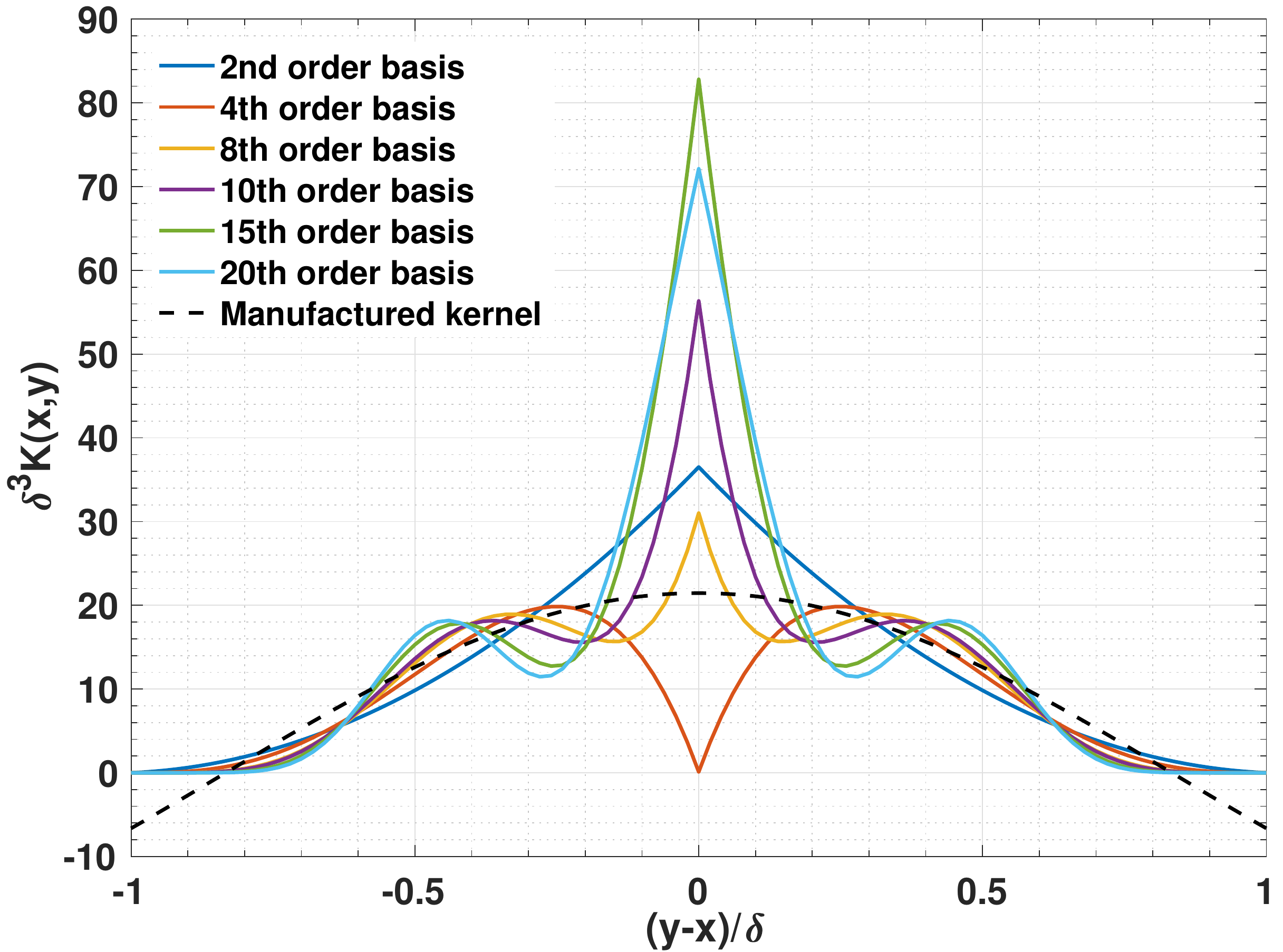}
    \includegraphics[width=0.48\textwidth]{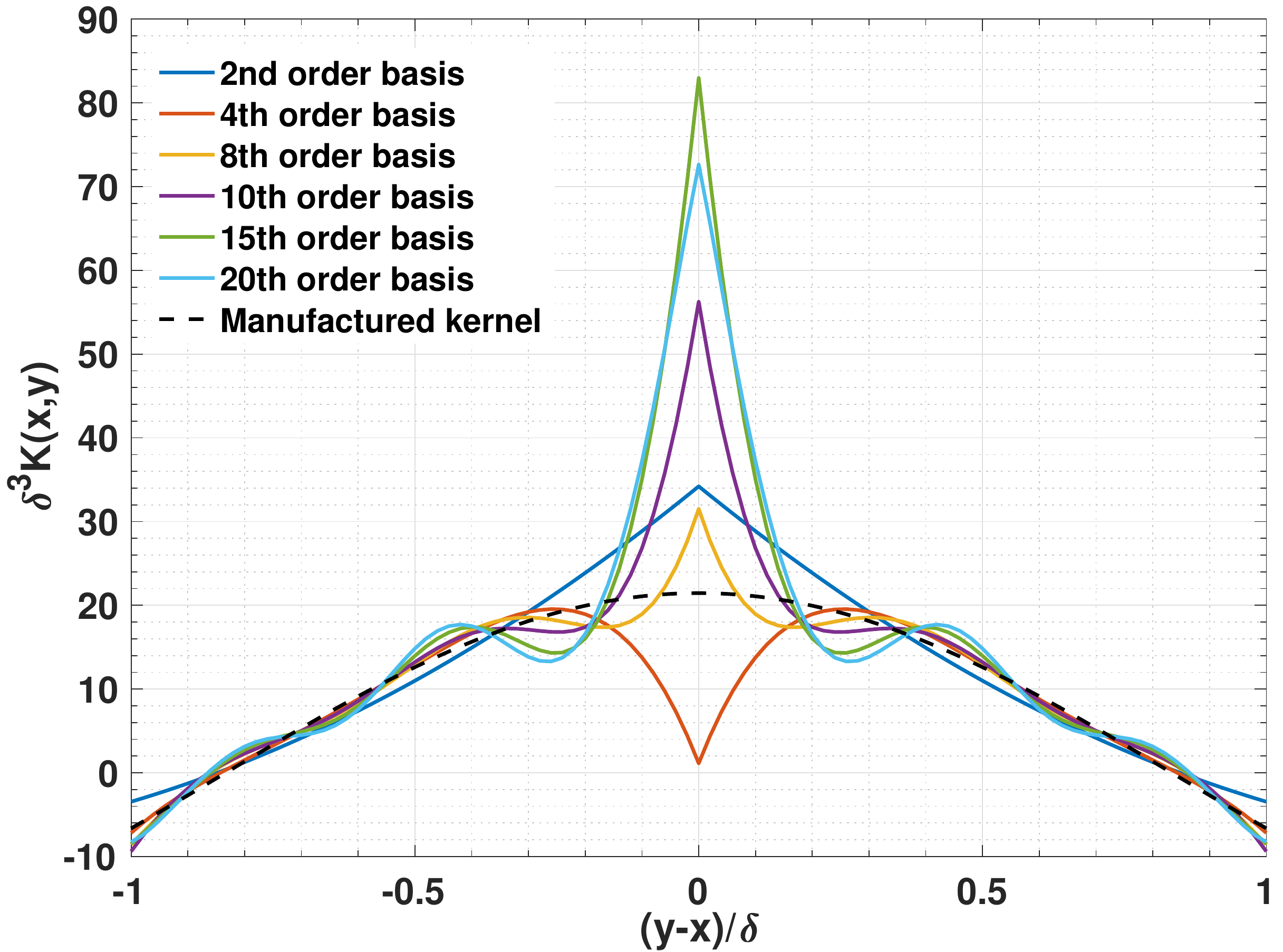}
    \caption{
\textit{Left: }Comparison of the nonnegative part $2\rho(x,y)$ of the kernel $K^*$ trained with low-frequency data and the manufactured kernel $K_{\text{man}}$ \eqref{eq:posnegkernel}. 
\textit{Right: }Comparison of the full kernel $K^* = 2\rho(x,y)+2h(x,y)$ trained with low-frequency data and the manufactured kernel $K_{\text{man}}$ \eqref{eq:posnegkernel}. Note the improved fit of the negative tails of the kernel.
}
    \label{fig:myreluC_i_and_CD_i}
\end{figure}

The plots of training and validation losses for both sets of data are shown in Figure \ref{fig:myrelulossvsorder_and_Loss_Hfreq}. This figure shows no benefit to increasing the basis order past $11$ for the case of low-frequency data, in contrast to the high-frequency test data for which both losses improve through order $20$. This test illustrates that the choice of basis order $M$ should take into account the frequency of the training and test data, an issue which will arise in Section \ref{sec:darcy} below. Next, to study the question of reproducing the manufactured kernel $K_{\text{man}}$ with the learned kernel $K^*$, we compare $K_{\text{man}}$ to $K^*$ first for the low-frequency data and varying basis order $M$ in Figure \ref{fig:myreluC_i_and_CD_i}. We  show both the nonnegative part $2\rho$ and the full kernel $2\rho + 2h$ to illustrate the contribution of $2h$.

\begin{figure}[htpb!]
    \centering
    \includegraphics[width=0.48\textwidth]{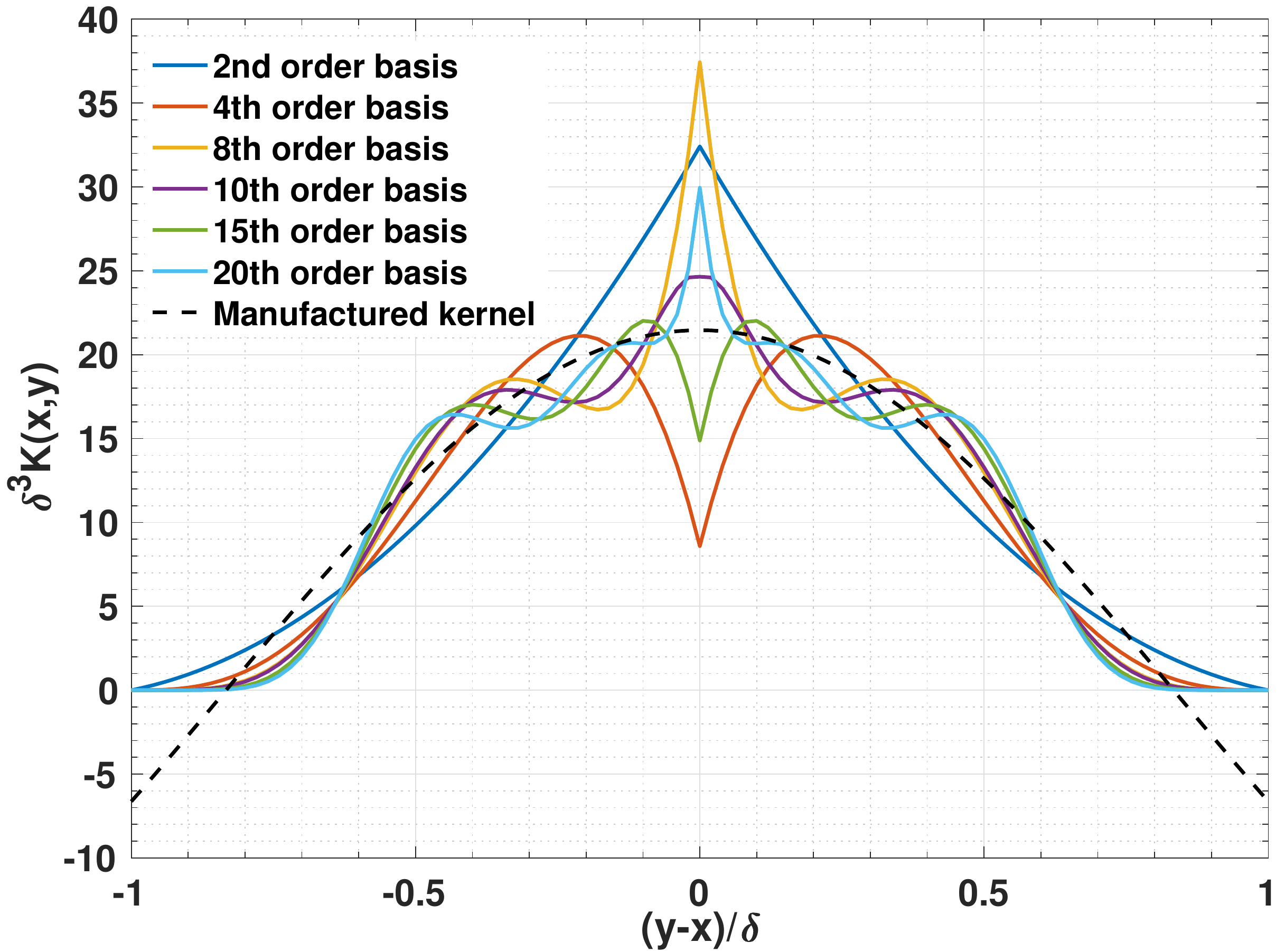}
    \includegraphics[width=0.48\textwidth]{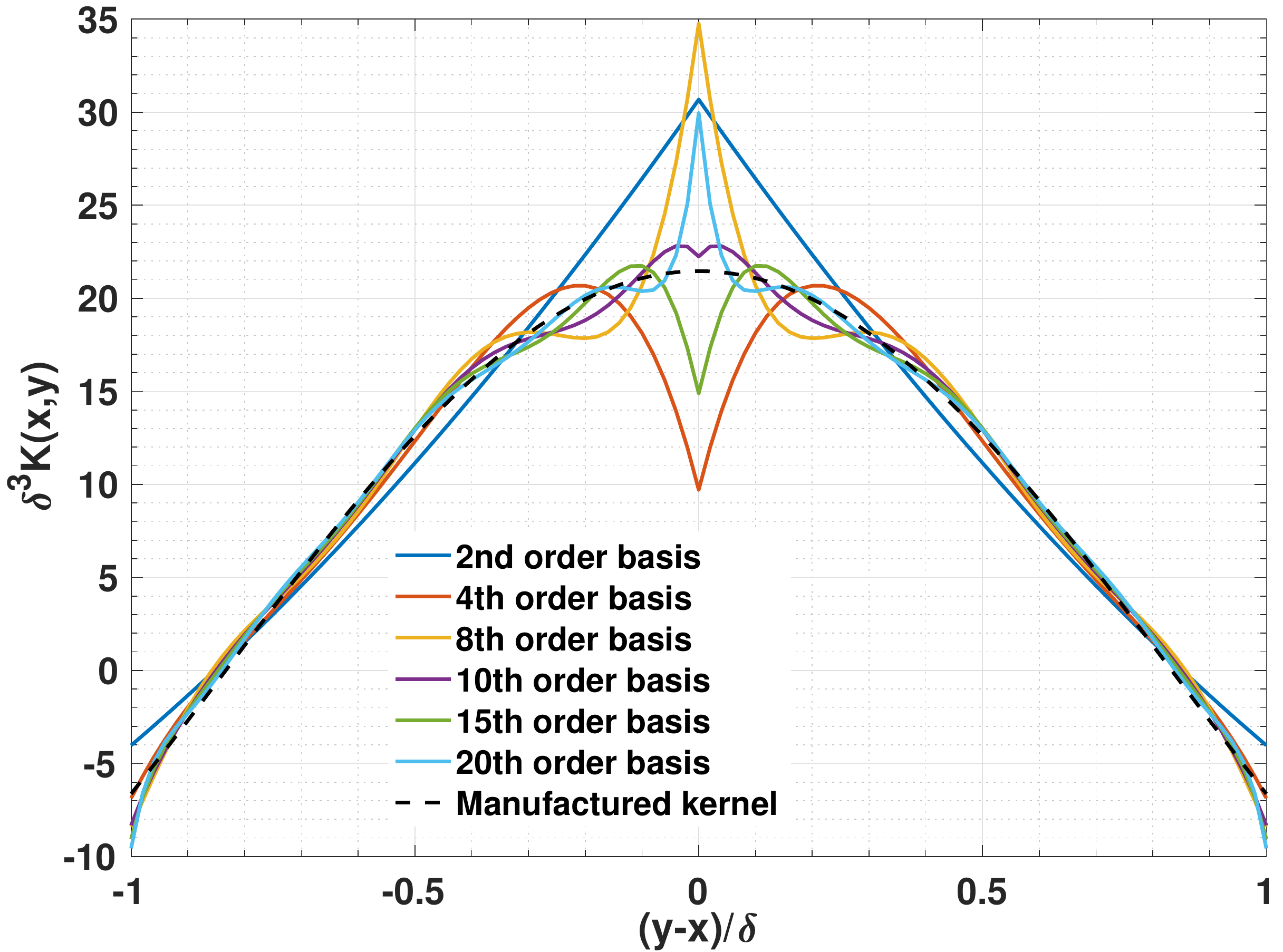}
    \caption{
\textit{Left: }Comparison of the nonnegative part $2\rho(x,y)$ of the kernel $K^*$ trained with high-frequency data and the manufactured kernel $K_{\text{man}}$ \eqref{eq:posnegkernel}. 
\textit{Right: }Comparison of the full kernel $K^* = 2\rho(x,y)+2h(x,y)$ trained with high-frequency data and the manufactured kernel $K_{\text{man}}$ \eqref{eq:posnegkernel}. Note the improved fit of the negative tails of the kernel.
}
    \label{fig:Hfreq_Ci_and_CDi}
\end{figure}

While the kernel is approximated by $2\rho + 2h$ for larger $|y-x|/\delta$ and sufficiently high-order basis, Figure \ref{fig:myreluC_i_and_CD_i} demonstrates that the manufactured kernel is not recovered for small $|x-y|/\delta$, implying that the action of the operator $\rho$ is insensitive to the kernel shape for small $|x-y|/\delta$ for training data $f_i$ of the form \eqref{eq:f_manufactured} with \eqref{eq:low_freq_distribution}. We hypothesize this is due to lack of higher frequencies in $f_i$ and the resulting solution $u_i$. By repeating this experiment in Figure \ref{fig:Hfreq_Ci_and_CDi} for the high-frequency data \eqref{eq:high_freq_distribution}, we obtain evidence for this hypothesis, as a closer fit is obtained for small $|x-y|/\delta$ than for the low-frequency training data in Figure \ref{fig:myreluC_i_and_CD_i} as the basis increases. 

These experiments suggest that our algorithm cannot be expected to recover exactly the kernel $K_{\text{man}}$ from the data $\mathcal{D}_{\text{train}}$. 
Moreover, resolution of the kernel in the vicinity of $|x-y|/\delta = 0$ is aided by the incorporation of higher frequency training data. Perhaps unsurprisingly, for such a sign-changing kernel, incorporation of $\bf D$ improves reproduction of $K_{\text{man}}$.

\section{Nonlocal coarse-graining of Darcy flow}\label{sec:darcy}

We next consider extraction of a nonlocal diffusion model by coarse-graining local solutions of Darcy's equation, a problem arising in subsurface flows through porous media \cite{darcy1856fontaines,whitaker1986flow}. In traditional inverse modeling and subsurface applications, the inference of an effective permissivity from data is a canonical problem motivating the development of multiscale and homogenization methods \cite{carrera2005inverse,payne2001effect,auriault2005transport,arbogast2012mixed}. Some works have argued that nonlocal models converging to Darcy's equations in the limit $\delta\rightarrow0$ provide an improved description of multiscale transport \cite{delgoshaie2015non,sen2012spatially}. We consider a simple 1D problem to demonstrate how nonlocal kernels may be derived from data, illustrating viability of data-driven multiscale nonlocal models for such processes. We numerically homogenize a repeated microstructure of lengthscale $2L$ using a nonlocal model with support lengthscale $\delta$. We will see that in this nonlocal context, the choice of coarse-graining lengthscale relative to $2L$ is tied to the effectiveness of the homogenization, similar to the local homogenization setting \cite{du2006size}.

\begin{figure}
\centering
\includegraphics[width=0.99\textwidth]{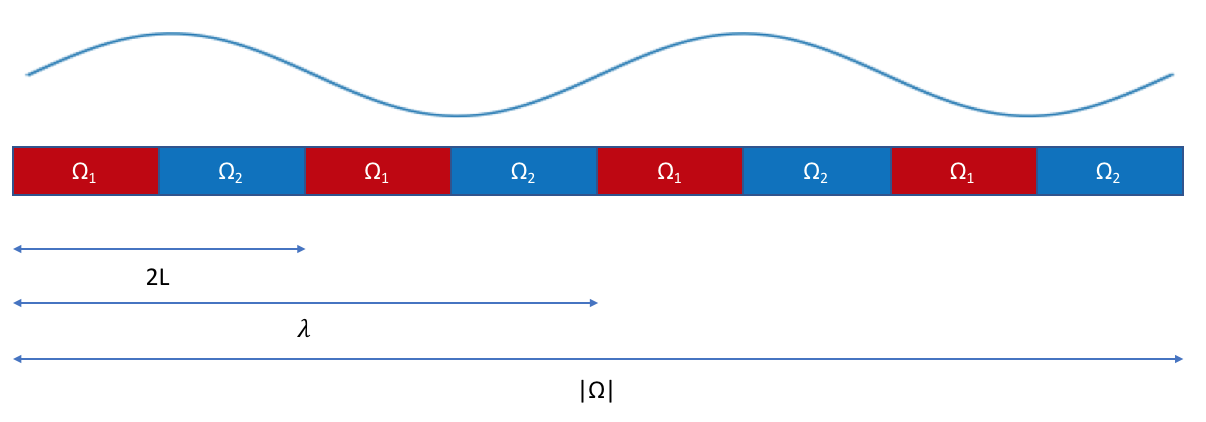}
\caption{Problem setup and relevant lengthscales for the Darcy coarse-graining problem. Periodic microstructure length is denoted $2L$, while forcing used during testing/training corresponds to $\sin(\lambda\,2 \pi x)$, and entire domain has measure $|\Omega|$. We will observe that choice of $\lambda$ used during test and training must be sufficiently large relative to $2L$ to obtain good fit for data-driven nonlocal model.}
\label{fig:darcySetup}
\end{figure}

We assume a periodic domain $\Omega = [0,10]$, constituted of alternating subdomains of width $L$ and piecewise constant diffusivity $\kappa_1 = 1$, $\kappa_2 = 4$ (Figure \ref{fig:darcySetup}). Further, we use a P1 nodal finite element solver to generate a collection of solutions $U = \left\{u_i\right\}_{i=1}^{N}$ to the following high-fidelity problem
\begin{align}\label{eqn:darcyTraining}
\nabla \cdot \mathbf{F}_i &= \sin(2 \pi x \lambda_i) := f_i\\
\mathbf{F}_i &= -\kappa(x) \nabla u_i.
\end{align}
Here, $\lambda_i$ denotes a wavelength sampled from the uniform distribution $\lambda_i \sim \mathcal{U}([\lambda_{\text{min}},{|\Omega|}])$. We omit details regarding the finite element solution, noting that we use $16$ elements per subdomain, which was chosen to ensure sufficient resolution of training data, and take $N=1000$ and $L=0.1$.

\begin{figure}
\centering
\includegraphics[width=0.49\textwidth]{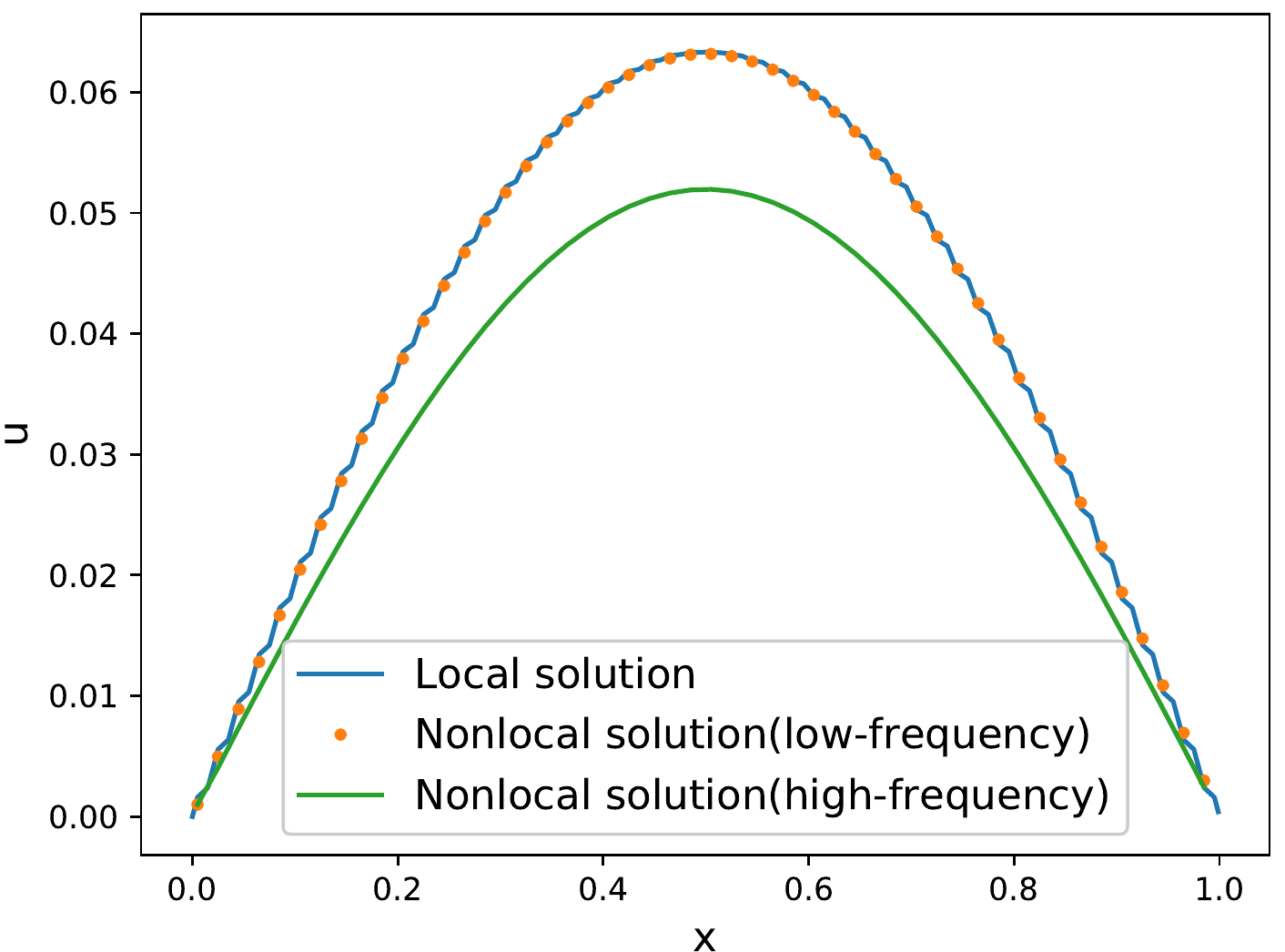}
\includegraphics[width=0.49\textwidth]{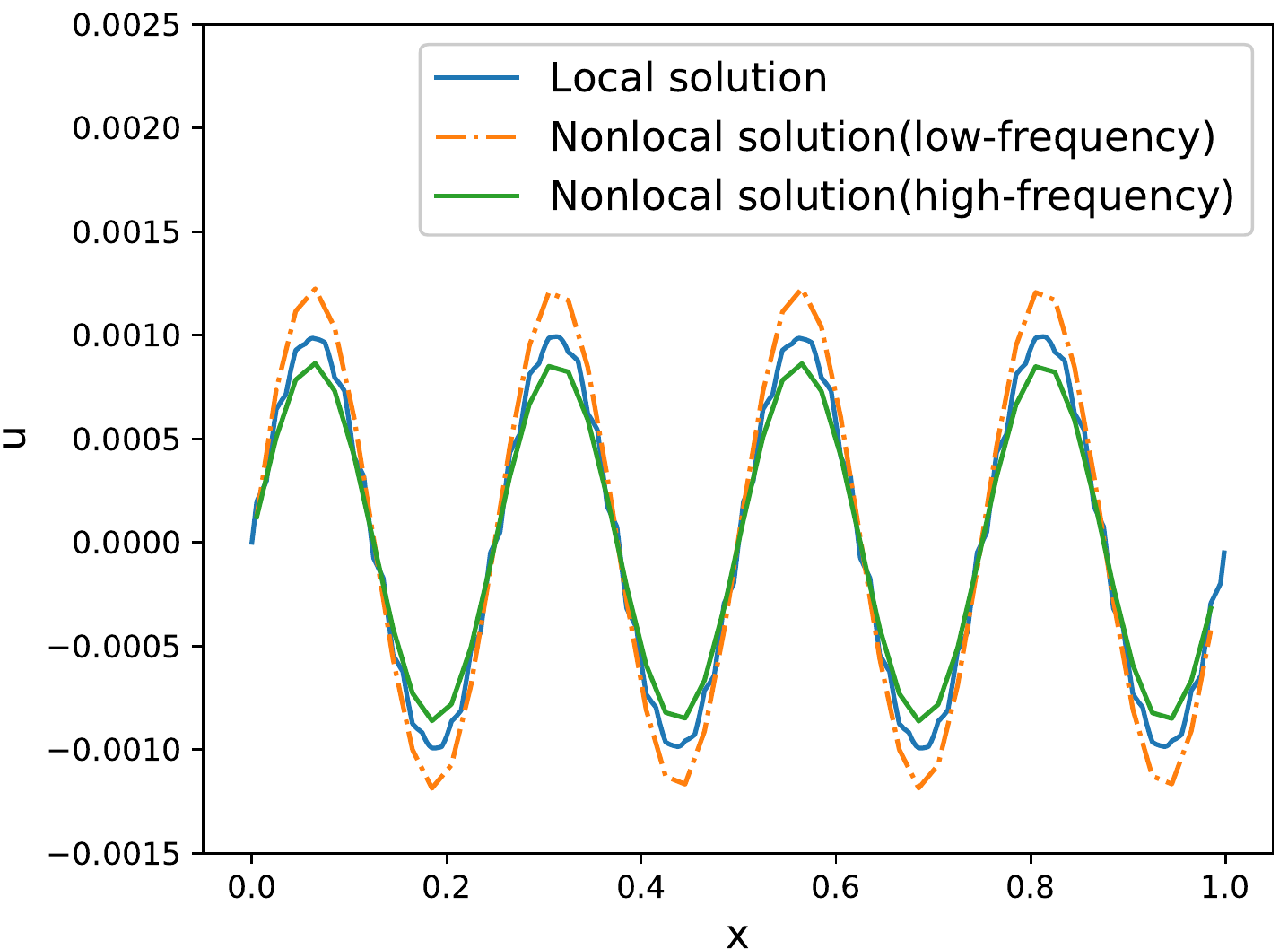}
\caption{Representative solutions to the extracted kernels, for $\delta = 8L$, testing against low-frequency \textit{(left)} and high-frequency \textit{(right)} forcing. We note that for low-frequency forcing, the solution is primarily smooth, and captured well by models extracted from low-frequency data. For high-frequency forcing, the local solution has artifacts from the microstructure, preventing extraction of an effective model using either high- or low-frequency training data.}
\label{fig:darcySolutions}
\end{figure}

We first coarse-grain $U$ into a collection of averaged solutions $\overline{U} = \left\{\overline{u}_i\right\}_{i=1}^{N}$, where we define $\overline{u}_i$ (and likewise $\overline{f}_i$) as piecewise constant over a given interval $x\in [ 2Ln,2L(n+1)]$,
\begin{equation}\label{eqn:coarsenedDarcy}
\overline{u}_i(x) = \frac{1}{2L}\int^{2L(n+1)}_{2Ln} u_i(y)\, dy,
\end{equation}
for a given domain $n$. In this manner, we obtain a representation of the local solution homogenized over the periodic microstructure. Finally, we apply Algorithm \ref{alg:augmented_lagrangian} to find a nonlocal kernel consistent with $\overline{U}$. That is, we apply the algorithm with 
\begin{equation}
\mathcal{D}_{\text{train}} = \{(\overline{u}_i, \overline{f}_i)\}_{i=1}^N
\end{equation}
to learn a nonlocal model \eqref{coarsegrained} with periodic boundary conditions.
The set $\mathcal{X}$ of test points in \eqref{eq:X_norm} is 500 uniformly distributed points from 0 to 10. 

A key feature of this data is the relative magnitude of the wavelength $\lambda_i$ and the microstructure $L$; similar to traditional homogenization we require $L/\lambda_i$ sufficiently small to treat the material as a homogeneous continuum. In Figure \ref{fig:darcySolutions} we present solutions to a pair of models with $\lambda_{\text{min}} = 1/16$ and $1$, both for $\delta = 8L$. When considering a low-frequency solution ($\mathcal{L}[u] = \sin \pi x$), the model trained on low-frequency data matches its local counterpart well, while for a high-frequency problem ($\mathcal{L}[u] = \sin 8 \pi x$) neither model performs well. This suggests that $L/\lambda_{\text{min}}$ must be sufficiently small in the training data to prevent the model from overfitting to high-frequency data. When applying the model, it will only agree with the local solution over sufficiently large lengthscales.

\begin{figure}[]
    \centering
    \includegraphics[width=.38\textwidth]{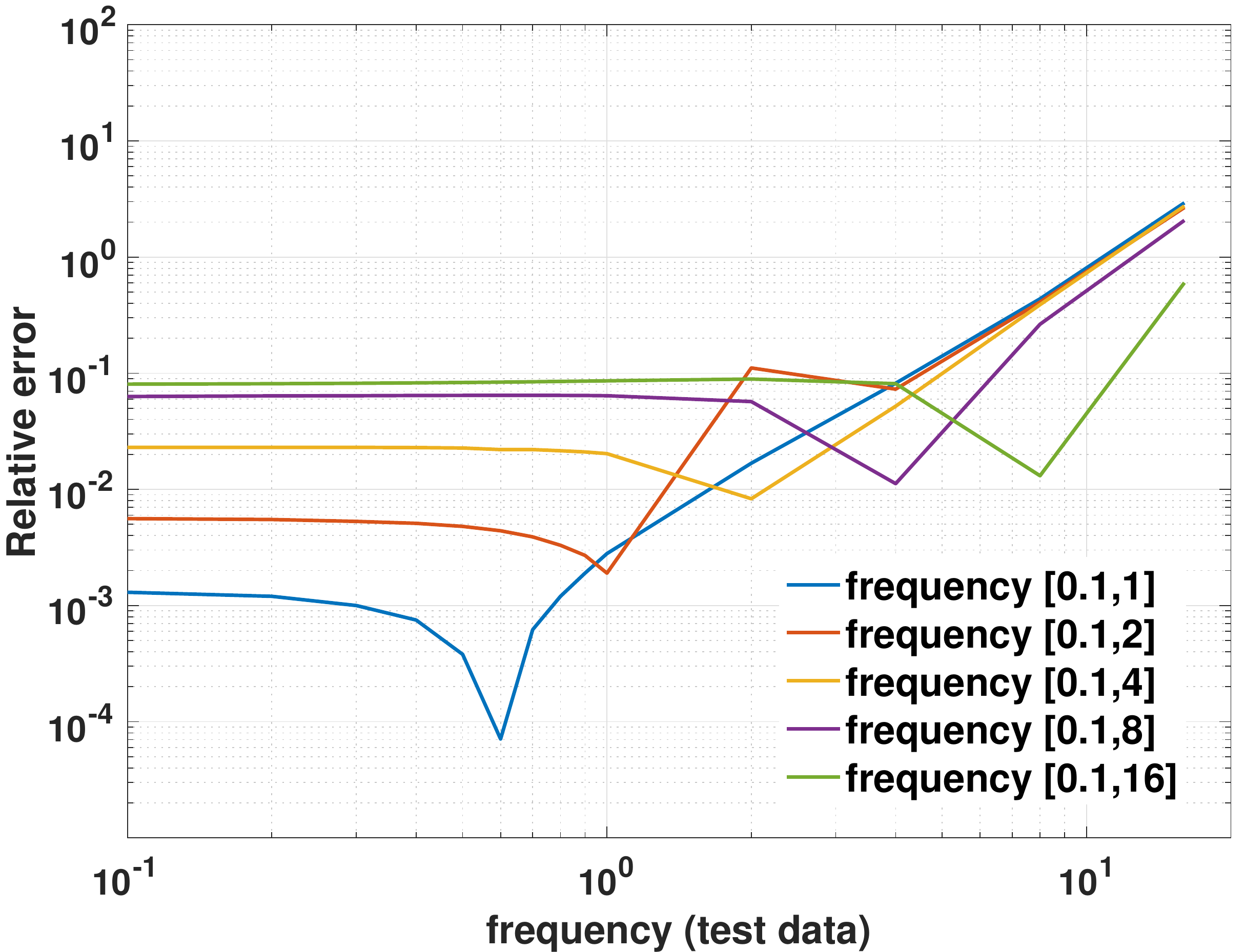}
    \includegraphics[width=.38\textwidth]{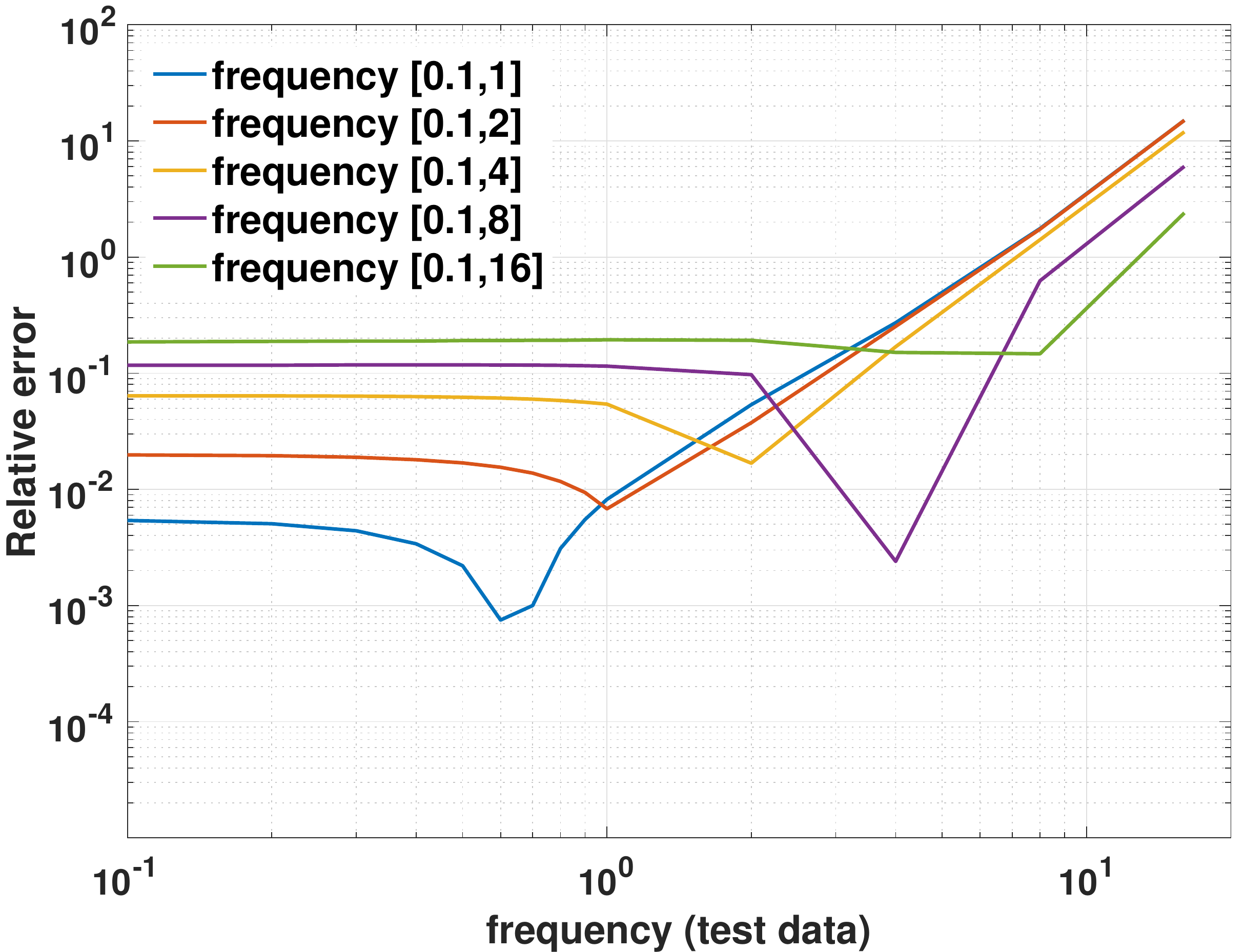}
    \includegraphics[width=.38\textwidth]{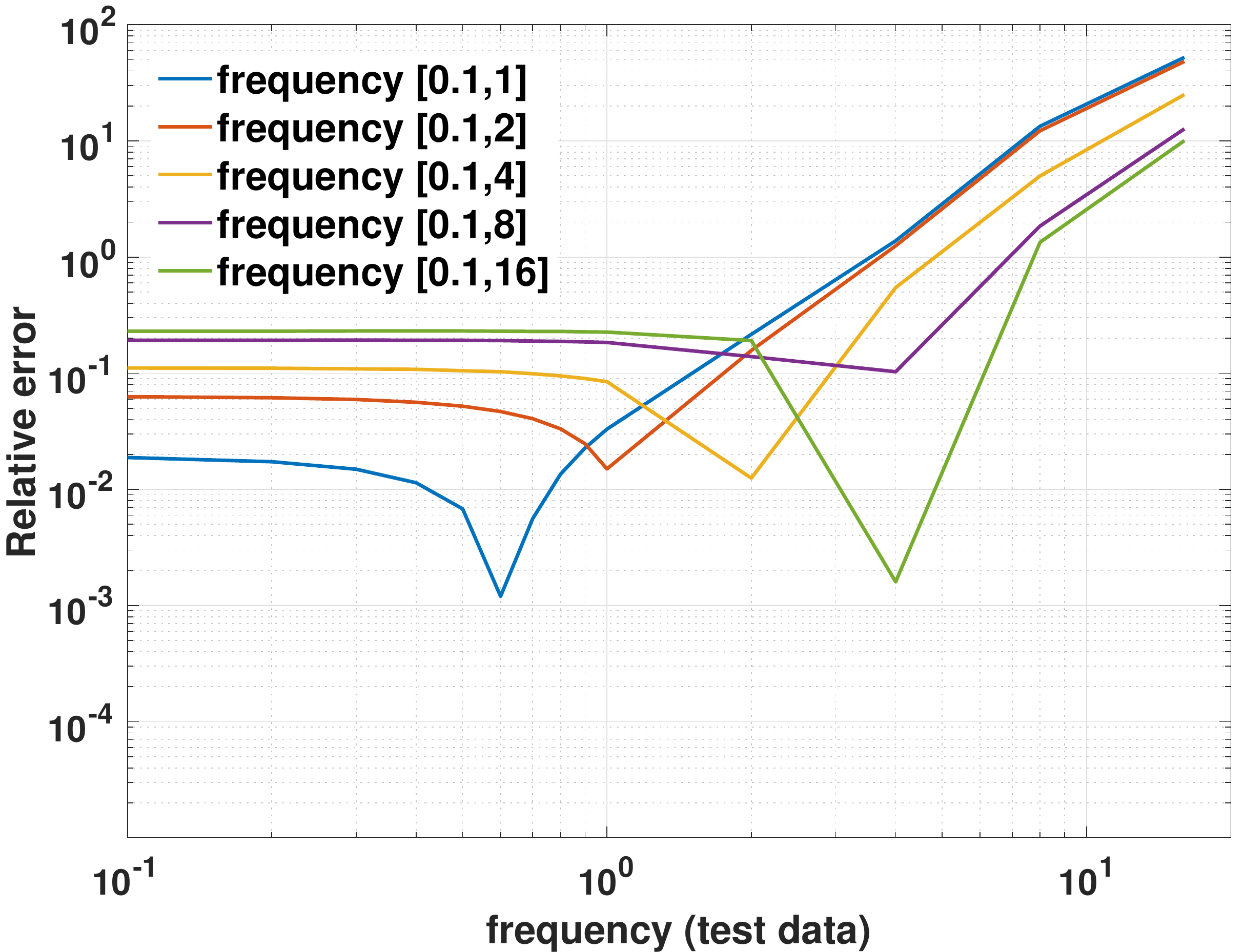}
    \vspace{-1ex}
    \caption{Relative solution error versus frequency of the forcing term in test data. The different curves in each subfigure signify a different range of frequency used in training data $[1/|\Omega|,1/\lambda_{\text{min}}]$. Increasing the ratio of $\lambda_{\text{min}}/L$ provides improved fit. \textit{Clockwise from top-left:} $\delta = 4L$, $\delta = 8L$, and $\delta = 16L$.}
    \label{fig:freq_sol16}
\end{figure}

\begin{figure}[]
    \centering
    \includegraphics[width=.38\textwidth]{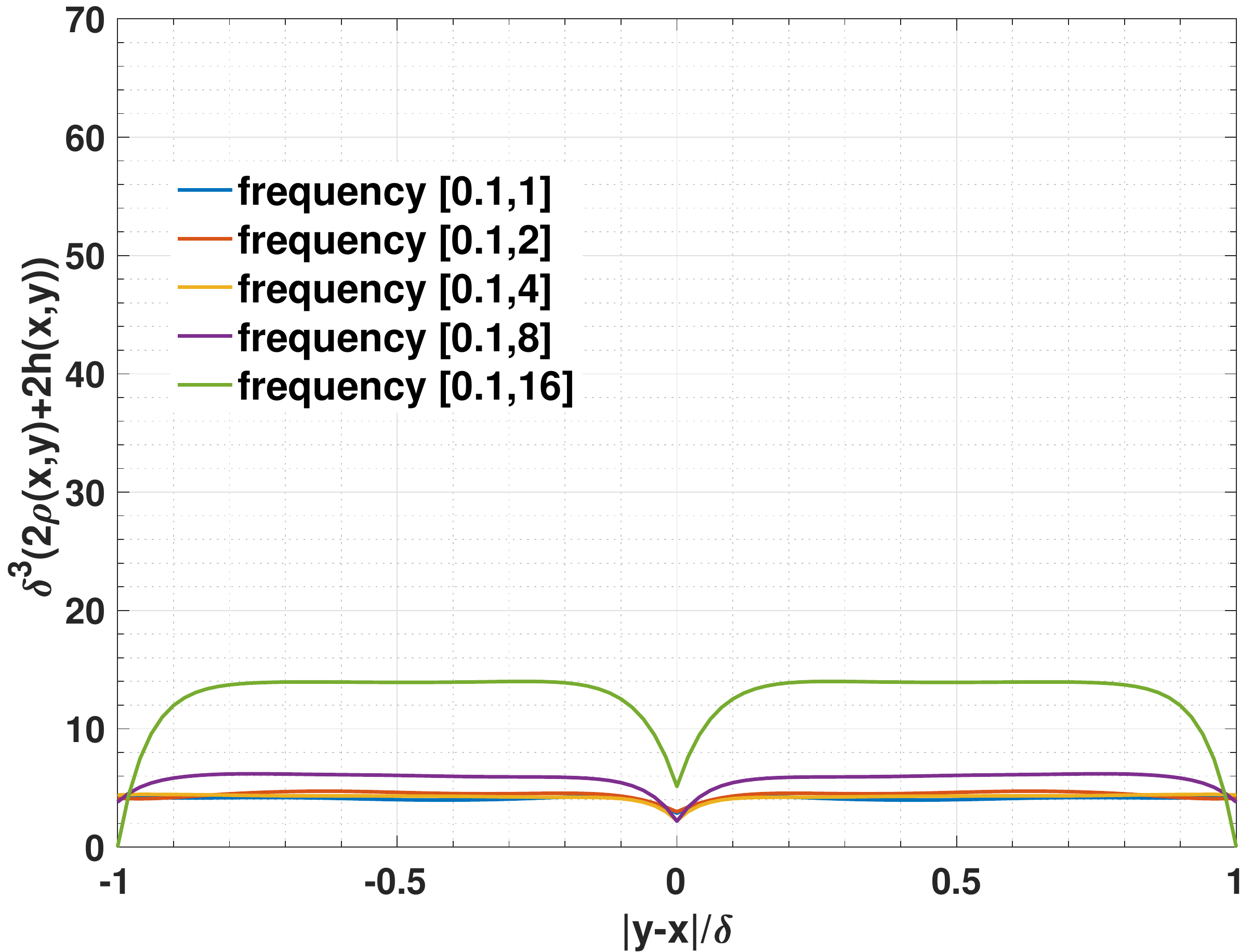}
    \includegraphics[width=.38\textwidth]{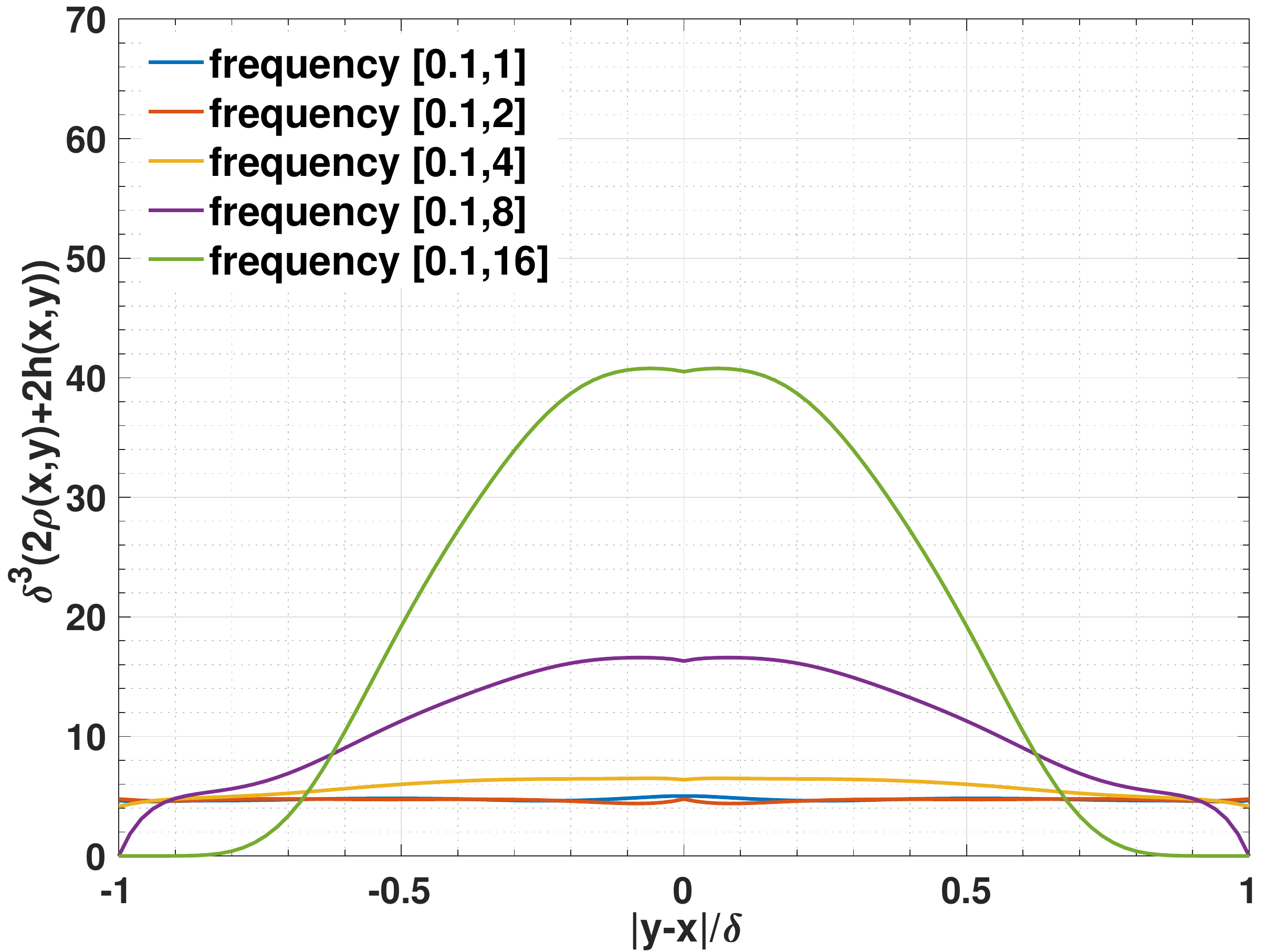}
    \includegraphics[width=.38\textwidth]{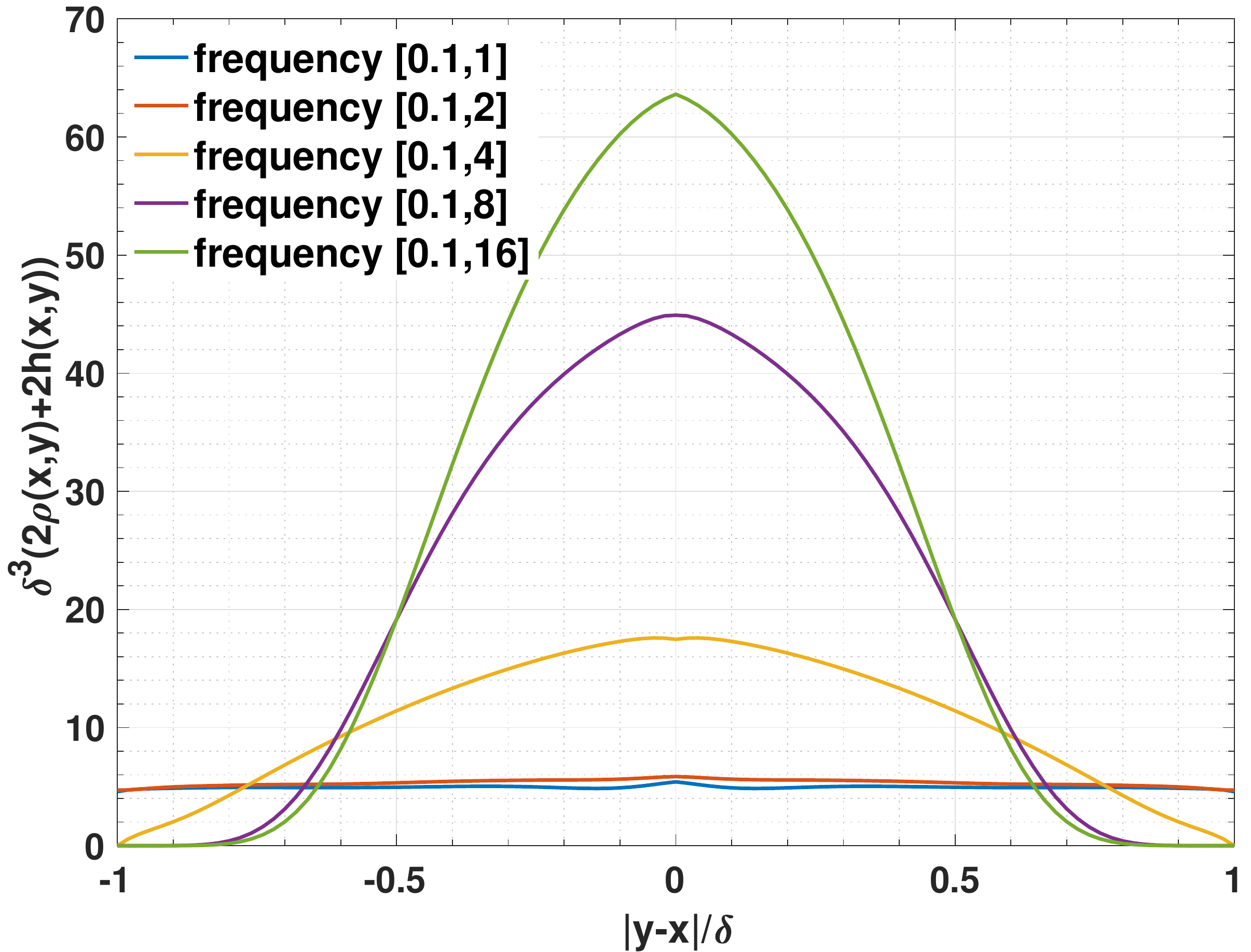}
    \vspace{-1ex}
    \caption{Trained kernels from different training data. Legend denotes range of frequency used in training data $[1/|\Omega|,1/\lambda_{\text{min}}]$. \textit{Clockwise from top-left:} $\delta = 4L$, $\delta = 8L$, and $\delta = 16L$ }
    \label{fig:darcykernels}
\end{figure}

Motivated by this observation, we perform a study sweeping over the choice of $\lambda_{\text{min}}$, both during training and testing, for $\delta \in \left\{4L,8L,16L\right\}$. The results of the study are presented in Figure \ref{fig:freq_sol16}, with the trend that smaller choice of $\delta$ provides generally better results. The extracted kernels associated with this study are presented in Figure \ref{fig:darcykernels}. We observe that for different $\delta$, the resulting kernel may have qualitatively different shape, and that for high-frequency training (which tests poorly) the kernel shifts to a different shape. Most notably, for this problem we observe no discernable negative part in the kernel, despite the fact that the algorithm allows for it. This suggests that, at least for this diffusion process, it is not necessary to introduce negativity into the kernel to obtain good fit;  similar nonnegative kernels were obtained in \cite{delgoshaie2015non} as an upscaling of a pore network.

\section{Extraction of nonlocal sign-changing kernels from high-order PDEs}\label{sec:homogenization}

We now seek an application where we expect to extract a sign-changing kernel. Weckner and Silling summarized several works where authors augment second-order elliptic problems with high-order derivatives to accommodate the high-frequency response of a given material \cite{weckner2011determination,eringen1972linear,eringen1992vistas,eringen1972nonlocal,jakata2008determination,mindlin1965second}. As an example, one may augment the description of a diffusion process by a Laplacian with higher-order operators such as biharmonics. Oftentimes these corrections provide a means to introduce lengthscales originating in microscales, e.g. Bazant introduces a correction of the form $c \delta^2 \nabla^4$ to the Poisson-Boltzmann equation to correct for the effect of finite ion size \cite{bazant2011double,storey2012effects}. Such local corrections may be related to nonlocal models by assuming a large degree of regularity, expanding the solution in a Taylor series, and matching the moments of the kernel to corresponding local terms. However, nonlocal models possess the desirable feature of only requiring $L^2$-regularity. Weckner and Silling derived in this manner a sign-changing kernel, where the high-frequency response manifests as a negative part. We will consider a simplified problem, and illustrate that the introduction of a very small negative component allows for a substantial improvement in model accuracy.

We seek a nonlocal surrogate for the biharmonic equation
\begin{align}\label{eq:biharmonic}
&u_{xx}-c\delta^2u_{xxxx} = f,\\
&u \text{ is periodic on } [0,L].
\end{align}
We generate data $\mathcal{D}_{\text{test}} = \{(u_i,f_i)\}_{i=1}^N$ for this equation for various $c$
and $\delta$ by randomly generating 
\begin{equation}\label{eq:biharmonic_f}
f_i = 
\sum_{k=1}^{99} \widehat{f}_k \cos(2\pi kx/L)
\end{equation}
for coefficients are sampled as
\begin{equation}\label{eq:biharmonic_distribution}
\widehat{f}_k \sim \exp(-\alpha k^2) \xi, \quad \xi \sim \mathcal{U}[0,1],
\end{equation}
where $\mathcal{U}[0,1]$ denotes the uniform distribution on $[0,1]$.
For each sampled $f_i$, the corresponding $u_i$ is solved from \eqref{eq:biharmonic} using the discrete Fourier transformation. We generate $N = 50,000$ pairs of $(u_i,f_i)$ for training, and train the kernel using an expansion \eqref{model-kernel} of order $M=20$ for a range of $\delta$ from $\delta = 0.125$ to $\delta = 0.99$ and $c$ from $c = 0.0001$ to $c = 0.1$. 
We train the kernel using $\mathcal{X}$ in \eqref{eq:X_norm} to be 101 equidistant points in $[0,1]$.

\begin{figure}[htpb!]
    \centering
    \includegraphics[width = .66\textwidth]{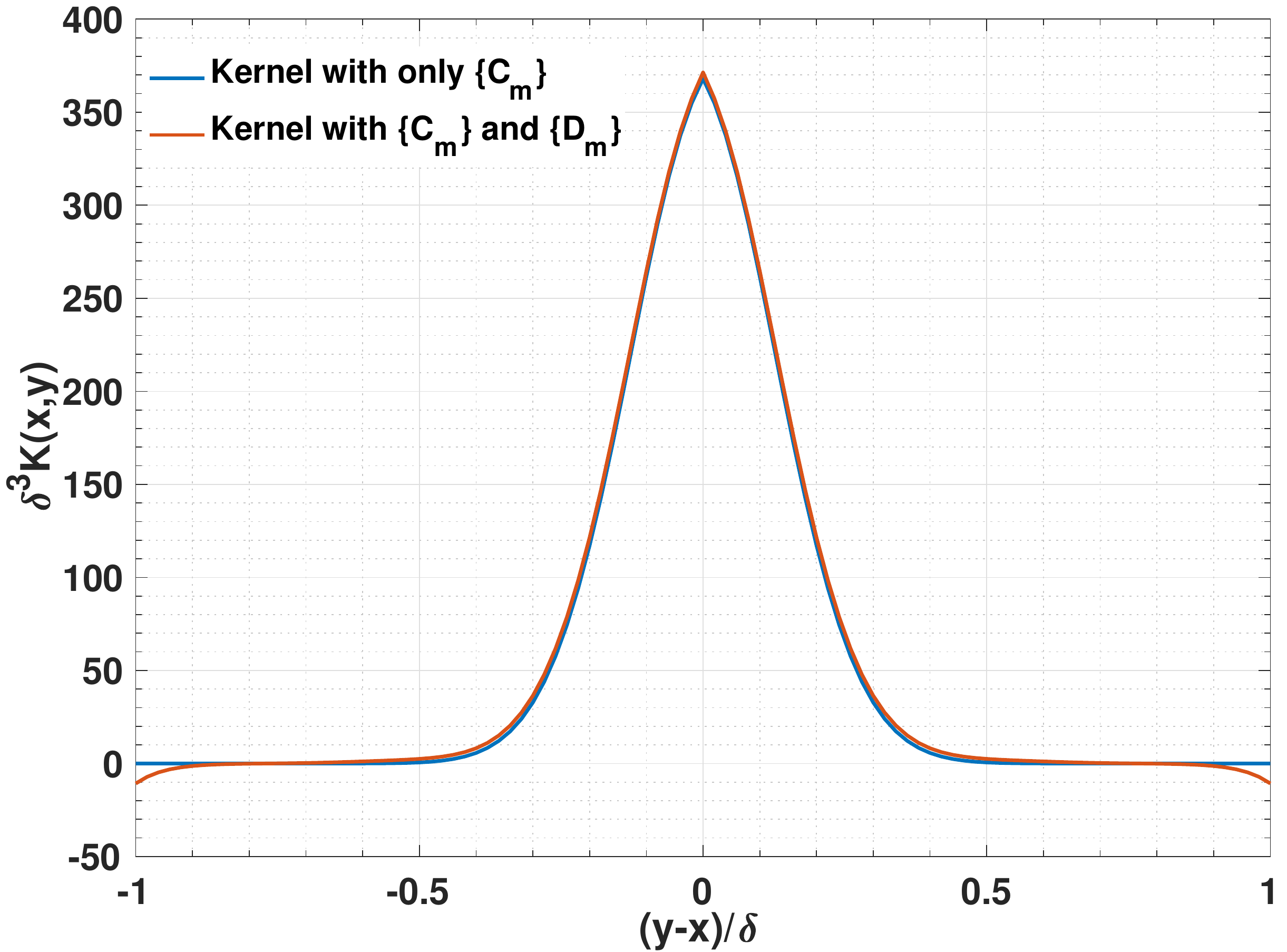}
    \caption{The comparison of the positive part of $K^*$ determined by $\{C_m\}$ and the full kernel $K^*$ given by both $\{C_m\}$ and $\{D_m\}$, when $c = 0.003$, $\delta =0.5$. Note that allowing non-positive kernels provides only a very small negative part near $|y-x| = \delta$, yet provides a substantial increase in accuracy (Tables \ref{tab:bihar_losses} and \ref{tab:bihar_solution}).}
    \label{fig:biharmonic_kernel}
\end{figure}

In Table \ref{tab:bihar_losses}, we report the losses after the first stage of the algorithm, when $\mathbf{C}$ has been trained and $\mathbf{D} = 0$, and the loss after the second stage of the algorithm when both $\mathbf{C}$ and $\mathbf{D}$ have been trained. Figure \ref{fig:biharmonic_kernel} compares these two kernels, for $c = 0.0003$ and $\delta = 0.5$, to illustrate that the two kernels are similar except for a small negative tail in the sign-changing kernel for large $|y-x|/\delta$. Nevertheless, Table \ref{tab:bihar_losses} illustrates that by fitting a sign-changing kernel as opposed to a nonnegative kernel, an average loss decreases by 19\% for $\delta = 0.125$, 44\% for $\delta = 0.25$, 42\% for $\delta = 0.5$, and 12\% for $\delta = 0.99$.

\begin{table}[htpb!]
\footnotesize
    \centering
    \begin{tabular}{|c|cc|cc|cc|cc|}
    \hline
        \multirow{2}{*}{$c$} & \multicolumn{2}{|c|}{$\delta = 0.125$}  & \multicolumn{2}{|c|}{$\delta = 0.25$}
        &\multicolumn{2}{|c|}{$\delta = 0.5$} & 
        \multicolumn{2}{|c|}{$\delta = 0.99$} \\
        \cline{2-9}
         & $\{C_m\}$ & $\{C_m + D_m\}$ &
          $\{C_m\}$ & $\{C_m + D_m\}$ & $\{C_m\}$ & $\{C_m + D_m\}$ &  $\{C_m\}$ & $\{C_m + D_m\}$ \\
         \hline
         0.0001 & 7.89e-5 & 6.49e-5 & 5.94e-4 & 1.48e-4 & 4.50e-3 & 1.50e-3 & 1.99e-2 & 1.78e-2 \\
         \hline
         0.0003 & 8.47e-5 & 8.45e-5 & 6.08e-4 & 1.46e-4 & 4.87e-3 & 1.88e-3 & 2.27e-2 & 1.71e-2 \\
         \hline 
         0.0005 & 9.46e-5 & 8.30e-5 & 6.79e-4 & 1.99e-4 & 5.58e-3 & 2.44e-3 & 2.46e-2 & 1.81e-2 \\
         \hline 
         0.001 & 9.97e-5 & 7.77e-5 & 7.17e-4 & 3.11e-4 & 6.12e-3 & 2.07e-3 & 2.88e-2 & 2.45-02 \\
         \hline 
         0.003 & 1.46e-4 & 1.32e-4 & 1.10e-3 & 5.14e-4 & 1.11e-2 & 4.76e-3 & 4.84e-2 & 3.90e-2 \\
         \hline 
         0.005 & 1.95e-4 & 1.74e-5 & 1.73e-3 & 1.32e-3 & 1.52e-2 & 9.16e-3 & 6.11e-2 & 5.46e-2 \\
         \hline 
         0.01 & 3.52e-4 & 3.14e-4 & 3.18e-3 & 1.94e-3 & 2.20e-2 & 1.64e-2 & 8.16e-2 & 7.20e-2 \\
         \hline 
         0.03 & 1.30e-3 & 1.25e-3 & 1.04e-2 & 8.09e-3 & 5.09e-2 & 4.32e-2 & 1.04e-1 & 1.08e-1 \\
         \hline 
         0.05 & 2.61e-3 & 2.38e-3 & 1.91e-2 & 1.51e-2 & 6.77e-2 & 5.48e-2 & 1.10e-1 & 1.18e-1 \\
         \hline 
         0.1 & 6.35e-3 & 6.18e-3 & 3.37e-2 & 3.47e-2 & 8.81e-2 & 7.83e-2 & 1.34e-1 & 1.23e-1 \\
         \hline 
    \end{tabular}
    \caption{Training loss when using $\{C_m\}$ only and $\{C_m+D_m\}$ under different settings. The loss averaged over all $c$ decreased by 19\%, 44\%, 42\%, and 12\% for $\delta = 0.125, 0.25, 0.5$ and $0.99$, respectively, when allowing for sign-changing kernels. }
    \label{tab:bihar_losses}
\end{table}

Table \ref{tab:bihar_solution} compares the relative error in solving \eqref{eq:biharmonic} for the test forcing
\begin{equation}\label{eq:test_case}
f = -4\pi^2\sin(2\pi x)-16c\delta^2\pi^4\sin(2\pi x)
\end{equation}
using the nonlocal model \eqref{model-kernel} with trained kernel $K^*$.
We solve the nonlocal equation $\mathcal{L}_{K^*} u = f$ with 
periodic boundary condition.
We compare the numerical solutions to the exact solution $u = \sin(2\pi x)$. 
This table corroborates Table \ref{tab:bihar_losses}, showing a significant in decrease testing error for $\delta = 0.25$ and $\delta = 0.5$ when using the sign-changing kernel rather than a strictly nonnegative kernel. 

\begin{table}[htpb!]
\footnotesize
    \centering
    \begin{tabular}{|c|cc|cc|cc|cc|}
    \hline
        \multirow{2}{*}{$c$} & \multicolumn{2}{|c|}{$\delta = 0.125$}  & \multicolumn{2}{|c|}{$\delta = 0.25$}
        &\multicolumn{2}{|c|}{$\delta = 0.5$} & 
        \multicolumn{2}{|c|}{$\delta = 0.99$} \\
        \cline{2-9}
         & $\{C_m\}$ & $\{C_m + D_m\}$ &
          $\{C_m\}$ & $\{C_m + D_m\}$ &  $\{C_m\}$ & $\{C_m + D_m\}$ & $\{C_m\}$ & $\{C_m + D_m\}$ \\
         \hline
         0.0001 & 0.49\% & 0.39\% & 1.31\% & 0.19\% & 3.15\% & 0.09\% & 5.74\% & 1.37\% \\
         \hline
         0.0003 & 0.48\% & 0.48\% & 1.27\% & 0.20\% & 3.39\% & 0.13\% & 5.66\% & 0.98\% \\
         \hline 
         0.0005 & 0.51\% & 0.54\% & 1.37\% & 0.25\% & 3.39\% & 0.13\% & 6.15\% & 1.42\% \\
         \hline 
         0.001 & 0.55\% & 0.44\% & 1.48\% & 0.46\% & 3.67\% & 0.12\% & 6.31\% & 0.87\% \\
         \hline 
         0.003 & 0.67\% & 0.61\% & 1.82\% & 0.52\% & 4.54\% & 0.15\% & 6.26\% & 0.86\% \\
         \hline 
         0.005 & 0.73\% & 0.71\% & 2.16\% & 1.11\% & 5.17\% & 0.24\% & 6.36\% & 0.94\% \\
         \hline 
         0.01 & 0.97\% & 1.04\% & 2.96\% & 1.05\% & 6.13\% & 0.42\% & 6.25\% & 1.27\% \\
         \hline 
         0.03 & 1.91\% & 1.67\% & 4.85\% & 2.59\% & 6.88\% & 0.20\% & 4.00\% & 0.97\% \\
         \hline 
         0.05 & 2.74\% & 2.33\% & 5.86\% & 2.97\% & 6.81\% & 0.19\% & 3.68\% & 0.71\% \\
         \hline 
         0.1 & 4.03\% & 3.69\% & 6.84\% & 5.17\% & 5.97\% & 1.29\% & 2.87\% & 0.95\% \\
         \hline 
    \end{tabular}
    \caption{Relative solution error for the test case \eqref{eq:test_case} when using $\{C_m\}$ only and $\{C_m+D_m\}$ under different settings. The test error averaged over all $c$ decreased by 9\%, 63\%, 94\%, and 80\% for $\delta =0.125, 0.25, 0.5$ and $0.99$, respectively, when allowing for sign-changing kernels. Note that while introduction of $D_m$ provides notable increase in accuracy, the resulting kernels differ by only a small negative tail (see Figure \ref{fig:biharmonic_kernel}).}
    \label{tab:bihar_solution}
\end{table}

\section{Fractional Laplacian}\label{sec:fractional}
In this section we consider a high-fidelity model substantially different from the ones considered in Sections \ref{sec:darcy} and \ref{sec:homogenization}: the Poisson problem for the fractional Laplacian \cite{Meerschaert2006}.
Fractional order equations are nonlocal equations often used in subsurface modeling; they accurately represent complex multi-scale phenomena by incorporating long-range interactions into the model itself in the form of a fractional-exponent derivative, as opposed to the classical integer-exponent derivative. This allows to capture the complete spectrum of diffusion through a few scalar model parameters. Despite being known to be the model of choice for scientific and engineering applications such as subsurface flow and transport \cite{Meerschaert2006}, their adoption is limited due to several technical challenges,  computational cost being the most important. This is understood by noting that the support of the nonlocal kernel is infinite, as opposed to the kernels considered in this paper whose support is limited to a ball of radius $\delta$. Formally, the fractional Laplacian operator is defined as 
\begin{equation}\label{eq:riesz_Laplacian_rn}
(-\Delta)^s u = C_{d,s}
\int_{\mathbb{R}^d} \frac{u(x) - u(y)}{|x-y|^{d+2s}}dy,
\end{equation}
where $s\in(0,1)$ and $C_{d,s}$ is a constant that depends on the dimension $d$ and the fractional order $s$. Here, the {\it positive} nonlocal kernel, with support on $\mathbb R^d$ is given by
\begin{equation}\label{eq:fractional-kernel}
K(x,y)= \dfrac{C_{d,s}}{|x-y|^{d+2s}}.
\end{equation}
We point out that, the interaction domain associated with the fractional Laplacian equation corresponds to the case of infinite horizon $\delta = \infty$.  Thus, the high-fidelity problem we consider is \eqref{high-fidelity} with $\mathcal{L}_{\text{HF}} = (-\Delta)^s$ and an exterior condition $u = q$ on $\Omega_I=\mathbb R^d\setminus\Omega$. In other words, exterior conditions must be prescribed on an infinite domain.

Our goal in this section is to answer the following question: can we learn a nonlocal compactly supported kernel (with a relatively small horizon) from high-fidelity data from a fractional model such that the action of the corresponding nonlocal operator is equivalent to the fractional Laplacian in \eqref{eq:riesz_Laplacian_rn}? A positive answer would enable accurate computation of fractional solutions at a much cheaper cost and increase the usability of fractional models now hindered by computational challenges.

To answer the question above, we first note that solutions to the Poisson problem for the \emph{truncated} fractional Laplacian (obtained by simply truncating the domain of integration in \eqref{eq:riesz_Laplacian_rn}) approach solutions of the non-truncated counterpart asymptotically with order $\mathcal O(\delta^{-2s})$ \cite{DElia2013}; this type of convergence holds in general for other instances of fractional operators such as the fractional gradient and divergence \cite{DElia2020}. Thus, for very small horizons, e.g., for $\delta\approx |\Omega|$, we cannot expect solutions to the truncated equation to be accurate approximations of solutions to the non-truncated one. This fact can be observed in the last column of Table \ref{tab:fractional} where we report the relative difference, for different values of $\delta$ between the non-truncated and truncated solutions; as expected, for very small $\delta$, the truncated solution is significantly different than the non-truncated one. These results are reported as a reference and their importance will be clear after we describe the impact of our algorithm. 

\smallskip
In order to reproduce the singularity in \eqref{eq:fractional-kernel} we slightly modify the kernel expansion \eqref{model-kernel} by dividing each basis function $B_{m,M}$ by $|x-y|^\alpha$ instead of $\delta^{d+2}$, where $\alpha$ is an additional tunable parameter. This means that we augment the set of parameters to be optimized within our algorithm by adding $\alpha$. Furthermore, due to the positivity of the fractional kernel \eqref{eq:fractional-kernel}, we only perform the first stage of Algorithm \ref{alg:augmented_lagrangian} to learn the coefficients $\bf C$ and the parameter $\alpha$. Recall that, as mentioned in the introduction, fractional equations are master equations of stochastic processes and, as such, their kernels represent (positive) jump rates.

We generate the training set $\mathcal D_{\rm train}=\{(u_i,f_i)\}_{i=1}^N$ by computing solutions of the fractional Laplacian equation with $s=0.75$ in $\Omega=B_1(0)$, a ball of radius $1$ and centered at $0$, and  exterior condition $q \equiv 0$ in $\mathbb{R} \setminus \Omega$. We sample $50,000$ forcing terms $f_i$ as in \eqref{eq:biharmonic_f} and \eqref{eq:biharmonic_distribution}, and obtain the corresponding solutions $u_i$ by evaluating the integral
\begin{equation*}
u_i(x) = \int_{B_1(0)} \mathcal G(x,y) f_i(y) dy,
\end{equation*}
using the trapezoidal quadrature rule,
where $\mathcal G$ is the Green's function \cite{DElia2013}.
We train with $\mathcal{X}$ in \eqref{eq:X_norm} consisting of 161 equidistant points in $[-0.8,0.8]$ to avoid singularities of the solutions $u_i$ near the boundary points $-1$ and $1$ \cite{Lischke2020}.
We test our algorithm for several values of $\delta$ and various order $M$ of the basis, initializing $\alpha$ to be $0$. In Table \ref{tab:fractional} we report the relative difference between the test solution for forcing $f=1$ and fractional order $s=0.75$ and the solution obtained by solving the nonlocal problem \eqref{coarsegrained} with kernel $K^*$. While a higher order does not necessarily yield a better prediction, for higher values of $\delta$, as expected, we obtain a much better approximation. As a reference we also report the difference that one would obtain by simply truncating the kernel \eqref{eq:fractional-kernel}. These results show that it is possible to obtain relatively accurate solutions to fractional equations at a much cheaper cost by using compactly supported fractional-type kernels learned from fractional data.
\begin{table}[htpb!]
    \centering
    \begin{tabular}{|c|c|c|c|c|c|}
    \hline
    \multirow{2}{*}{$\delta$} &  \multicolumn{5}{c|}{relative difference} \\
    \cline{2-6}
     & order = 0 &   order = 5 & order = 10 & order =20 & truncated kernel\\
    \hline
        $2^{-3}$ & 9.46\% &7.93\% &16.8\% & 19.8\% & 319\%\\
        $2^{-2}$ & 15.4\%& 22.89\% & 23.5\%& 26.0\% & 182\%\\
        $2^{-1}$ & 13.2\%& 8.57\% & 10.0 \% &7.2\% & 98\%\\
        $2^{0}$  & 3.97\%& 1.76\%& 2.37\%& 2.88\% & 48\%\\
        $2^{1}$ & 0.84\%& 2.18\%& 1.84\%&3.43\% &23\%\\
        $2^{2}$ & 1.59\%& 4.25\%& 3.66\%&1.95\% & 15\%\\
        \hline
    \end{tabular}
    \caption{Relative difference between the analytic solution with forcing $f=1$ and fractional order $s=0.75$ and the solution obtained in correspondence of the reconstructed kernel for several values of $\delta$ and different basis orders.}
    \label{tab:fractional}
\end{table}

\section{Conclusion}\label{sec:conclusion}

In this work we have presented an optimization framework for discovering sign-changing nonlocal models from high-fidelity synthetic data. In the nonlocal community there are several open questions about the derivation of such models from first-principles, and the data-driven approach presented here can provide guidance regarding what models emerge naturally from data. The fundamental property we have pursued is a guarantee that learned models be solvable - this ensures robustness without requiring access to complicated PDE-constrained optimization codes, and may be implemented in popular optimization packages such as Tensorflow \cite{tensorflow2015-whitepaper} and PyTorch \cite{NEURIPS2019_9015}. While we worked with Bernstein polynomials for ease of implementing the resulting inequality constraints, an interesting area of future research would be to consider whether incorporating deep learning architectures into this framework provides benefit. Working with polynomials assumes regularity in the underlying kernels, while e.g. shallow ReLU networks would allow parameterization of discontinuous kernels \cite{cyr2019robust,he2018relu}.

While we have focused on simple one-dimensional experiments for ease of presentation and to ensure examples are easily reproducible, the results demonstrate the potential impact this framework may have for range of important problems. We have demonstrated in Section \ref{sec:darcy} how one may use this approach to perform coarse-graining without incorporating restrictive and mathematically complex derivations. In Section \ref{sec:homogenization}, we illustrated the instrumental role sign-changing kernels have deriving reduced-regularity nonlocal versions of high-order PDEs, which have been shown previously to be fundamental in resolving high-frequency response of certain solid materials \cite{weckner2011determination}. In Section \ref{sec:fractional}, we see that one may extract computationally efficient compactly supported models providing sparse discrete operators that accurately approximate fractional operators with infinite horizons. This is particularly promising as a means of deriving preconditioners and $O(n)$ solvers for fractional systems, as naive discretizations provide dense matrices which require complex hierarchical preconditioners to solve efficiently \cite{ainsworth2018towards,xu2018efficient,Glusa2019}. However, these examples provide only a first sample of possible directions to learn data-driven nonlocal models, and we pursue in future work application of this framework toward large-scale datasets more representative of open problems in science and engineering. 

\section*{Acknowledgements}
Sandia National Laboratories is a multi-mission laboratory managed and operated by National Technology and Engineering Solutions of Sandia, LLC., a wholly owned subsidiary of Honeywell International, Inc., for the U.S. Department of Energy’s National Nuclear Security Administration under contract {DE-NA0003525}. This paper describes objective technical results and analysis. Any subjective views or opinions that might be expressed in the paper do not necessarily represent the views of the U.S. Department of Energy or the United States Government. Report Number: SAND2020-5158 J.

The work of N. Trask, M. Gulian, and M. D'Elia is supported by the U.S. Department of Energy, Office of Advanced Scientific Computing Research under the Collaboratory on Mathematics and Physics-Informed Learning Machines for Multiscale and Multiphysics Problems (PhILMs) project, and by the Sandia National Laboratories Laboratory Directed Research and Development (LDRD) program. M. Gulian is supported by the John von Neumann fellowship at Sandia National Laboratories. H. You and Y. Yu are supported by the National Science Foundation under award DMS 1753031.

\bibliographystyle{elsarticle-num}
\bibliography{snl}

\end{document}

%% file: mathhead.tex
\usepackage{amsmath}
\usepackage{amsfonts}
\usepackage{amssymb}
\usepackage{amsthm}
\usepackage{bm}
\usepackage{indentfirst}
\usepackage{graphicx}
\usepackage{subfigure}
\usepackage{array}
\usepackage{fullpage}
\usepackage{color}

\DeclareMathAlphabet{\mathsfsl}{OT1}{cmss}{m}{sl}

\newcommand{\PreserveBackslash}[1]{\let\temp=\\#1\let\\=\temp}
\newcolumntype{C}[1]{>{\PreserveBackslash\centering}p{#1}}
\newcolumntype{R}[1]{>{\PreserveBackslash\raggedleft}p{#1}}
\newcolumntype{L}[1]{>{\PreserveBackslash\raggedright}p{#1}}

\renewcommand{\vartheta}{\Theta}

\definecolor{mygreen}{rgb}{0.1,0.75,0.2}

\numberwithin{equation}{section}

\theoremstyle{definition}